\documentclass{article}
\pdfoutput=1
\usepackage{amsmath,amssymb,amsthm,picins,color,graphicx}
\usepackage{hyperref}
\usepackage[margin=1.3in]{geometry}
\hypersetup{
pdfauthor = {Boris Bukh},
pdftitle = {Measurable sets with excluded distances},
pdfsubject = {Mathematics},
pdfkeywords = {excluded distances, Hadwiger-Nelson problem,
chromatic number of the plane, measurable coloring, distance graph,
chromatic number, Euclidean Ramsey theory},
pdfstartview = {FitH},
pdfpagemode = {None}}

\newcounter{saveenum}

\newcommand*{\R}{\mathbb{R}}
\newcommand*{\Q}{\mathbb{Q}}
\newcommand*{\Z}{\mathbb{Z}}

\newcommand*{\abs}[1]{\lvert #1\rvert}
\newcommand*{\meas}{\abs}

\newcommand*{\norm}[1]{\lVert #1\rVert}
\newcommand*{\floor}[1]{\lfloor #1\rfloor}
\newcommand*{\Zm}{\mathcal{Z}}
\newcommand*{\Msr}{\mathfrak{M}}
\newcommand*{\proofpart}[1]{\textbf{#1: }}
\newtheorem{theorem}{Theorem}
\newtheorem{lemma}[theorem]{Lemma}

\newtheorem{corollary}[theorem]{Corollary}
\newtheorem{definition}[theorem]{Definition}

\newtheorem*{saturationlemma}{Supersaturation lemma}
\newtheorem*{saturationlemmaalt}{Weak supersaturation lemma}
\DeclareMathOperator{\diam}{diam}
\DeclareMathOperator{\supp}{supp}
\DeclareMathOperator{\dist}{dist}
\newcommand*{\Lem}{\mathtt{Lemma}}
\newcommand*{\WkLem}{\mathtt{WeakLemma}}
\newcommand*{\Lembf}{\mathtt{\mathbf{Lemma}}}
\newcommand*{\WkLembf}{\mathtt{\mathbf{WeakLemma}}}
\newcommand*{\OR}{\text{ OR }}
\newcommand*{\AND}{\text{ AND }}

\renewcommand*{\epsilon}{\varepsilon}

\author{Boris Bukh}
\title{Measurable sets with excluded distances\protect\footnote{This work
is a part of a Ph.~D.\ thesis under the supervision of Benjamin Sudakov.}}

\begin{document}
\maketitle
\begin{abstract}
For a set of distances $D=\{d_1,\dotsc,d_k\}$ a set
$A$ is called $D$-avoiding if no pair of points
of $A$ is at distance $d_i$ for some~$i$. We show that
the density of $A$ is exponentially small in $k$
provided the ratios $d_1/d_2$, $d_2/d_3$, \dots, $d_{k-1}/d_k$
are all small enough. This resolves a question
of Sz\'ekely, and generalizes a theorem
of Furstenberg-Katznelson-Weiss, Falconer-Marstrand, and Bourgain.
Several more results on $D$-avoiding sets are presented.
\end{abstract}

\section{Introduction}
The problem of determining the least number of colors required to
color the points of the plane $\R^d$ so that no pair of points
at distance $1$ is colored in the same color was first investigated by
Nelson and Hadwiger in 1940s. This number, which
we denote by $\chi_{\R^d}(\{1\})$, is called the chromatic
number of $\R^d$ because it is the chromatic number of the graph
whose vertices are the points of $\R^d$ and the edges are
pairs of points that are distance $1$ apart. We denote this
graph by $G_{\R^d}(\{1\})$.

In the dimension two, there has been no improvement on the 
bounds $4\leq \chi_{\R^2}(\{1\})\leq 7$ in the past
forty-five years \cite{cite:hadwiger,cite:moser_moser}.
In higher dimensions, however, Frankl and Wilson
\cite{cite:frankl_wilson} showed that the chromatic number
grows exponentially in the dimension, 
$\chi_{\R^d}\bigl(\{1\}\bigr)\geq (1.207\ldots+o(1)\bigr)^d$, 
confirming an earlier conjecture of Erd\H{o}s. The paper
of Frankl and Wilson in conjunction with the earlier
work of Ray-Chaudhuri and Wilson\cite{cite:rayc_wilson} 
laid down the theory of set families with restricted
intersection, which led to many other results including 
the disproof of Borsuk's conjecture by Kahn and Kalai \cite{cite:kahn_kalai}.

It was first shown by Erd\H{o}s and de Bruijn \cite{cite:erdos_debruijn}
that the chromatic number of any infinite graph, and $G_{\R^d}(\{1\})$
in particular, is the maximum of the chromatic numbers of its finite
subgraphs, provided the maximum is finite. The proof relied on the axiom
of choice, which suggested that the chromatic number might depend
on the underlying axiom system. This was partially confirmed by
Falconer \cite{cite:falconer_coloring} who showed that
there is no coloring of $\R^2$ into four colors such that
each color class is a Lebesgue measurable set and
no pairs of points at distance $1$ have the same color. Since as
shown by Solovay \cite{cite:solovay} the axiom that all
subsets of $\R$ are Lebesgue measurable is consistent with 
the usual Zermelo-Fraenkel set theory without the axiom of choice, 
$\chi_{\R^2}(\{1\})=4$ is unprovable in the set theory
without the axiom of choice.

Thus, we denote by $\chi^m_{\R^d}(\{1\})$ the least number of colors
required to color $\R^d$ so that no points at distance $1$ are assigned
the same color, and each color class is a measurable set. A set
with no pairs of points at distance $1$ is going
to be called \emph{$\{1\}$-avoiding}.
The most natural way to show that $\chi^m_{\R^d}({\{1\}})$ is large is 
by showing that no color class can be large. 
Denote by $\bar{d}(A)$ the upper limit density of $A$ (which is
formally defined in section \ref{sec_notation}). Let 
$m_{\R^d}(\{1\})=\sup \bar{d}(A)$ be the supremum
over all measurable $\{1\}$-avoiding sets. Then 
$\chi^m_{\R^d}(\{1\})\leq 1/m_{\R^d}(\{1\})$. Unfortunately,
Falconer's proof that $\chi^m_{\R^2}(\{1\})\geq 5$ does not show
that $m_{\R^2}(\{1\})<1/4$. 
The best known bounds are $0.229365\leq m_{\R^2}(\{1\})\leq
12/43$ (see \cite[p.~61]{cite:scheinerman_ullman_book} and 
\cite{cite:szekely_thesis} respectively), and it is a 
conjecture of Erd\H{o}s that $m_{\R^2}(\{1\})<1/4$ \cite{cite:szekely}.

The problem of forbidding more than one distance was first
studied by Sz{\'e}kely in his thesis \cite{cite:szekely_thesis}. There
he established the first bounds on $\chi^m_{\R^d}(D)$ and $m_{\R^d}(D)$
which denote the analogues of $\chi^m_{\R^d}(\{1\})$ and $m_{\R^d}(\{1\})$,
respectively, where a finite set of distances $D=\{d_1,\dotsc,d_k\}$
is forbidden. Sz\'ekely conjectured that in dimension $d\geq 2$ for 
any set $A$ with $\bar{d}(A)>0$ there is a $d_0$ such that all the distances
greater than $d_0$ occur among the points of~$A$. The conjecture
was proved by Furstenberg, Katznelson and Weiss \cite{cite:fkw}. Their proof
was ergodic-theoretic. Later Bourgain found a harmonic-analytic proof
\cite{cite:bourgain}, and Falconer and Marstrand gave a direct
geometric proof \cite{cite:falconer_marstrand}. Sz\'ekely also conjectured 
that if $d_1,d_2,\dotsc$ is a sequence converging to $0$, then 
$m_{\R^d}(\{d_1,\dotsc,d_k\})\to 0$ as $k\to\infty$. This was proved
by Falconer \cite{cite:falconer_small} and Bourgain \cite{cite:bourgain}.

It is not known how large $\chi_{\R^d}(D)$ can be for a set
$D$ of given size. It has been long known 
that $\sup_{\abs{D}=k} \chi_{\R^2}(D)\geq c k\sqrt{\log k}$ 
\cite[p.~180]{cite:cfg_unsolved}. The only known upper 
bound $\sup_{\abs{D}=k} \chi_{\R^d}(D)\leq \chi_{\R^d}(\{1\})^k$ comes from
the observation that the coloring, which is a product of colorings that avoid
$D_1$ and $D_2$, avoids both $D_1$ and $D_2$. Croft, Falconer and Guy asked
whether $\sup_{\abs{D}=k}\chi_{\R^d}(D)$ is exponential in $k$ \cite[Prob.~G11]{cite:cfg_unsolved}. Erd\H{o}s conjectured that 
$\sup_{\abs{D}=k}\chi_{\R^d}(D)$ is polynomial in $k$ 
\cite{cite:erdos_polygrowth}.

In this paper we answer the question of Croft, Falconer and Guy in
the measurable setting by showing that in the dimension $d\geq 2$ 
as the ratios $d_1/d_2,d_2/d_3,\dotsc,d_{k-1}/d_k$ all tend to infinity 
$m_{\R^d}(D)$ tends to $m(\{1\})^k$, and thus $\sup_{\abs{D}=k}\chi^m_{\R^d}(D)
\geq 1/m(\{1\})^k$. We will also
show that $m_{\R^d}(D)\geq m(\{1\})^k$ for every set of $k$ distances $D$,
answering question of Sz\'ekely 
\cite[p.~657]{cite:szekely}, who asked for the value of $\inf_{\abs{D}=k} m_{\R^d}(D)$.
This also generalizes the above-mentioned 
theorems of Furstenberg-Katznelson-Weiss and Falconer. Indeed,
to deduce Furstenberg-Katznelson-Weiss theorem suppose there is
a set $A$ with $\bar{d}(A)>0$ and a sequence $d_1,d_2,\dotsc$ going
to infinity such that the distance $d_i$ does not occur between
points of $A$. Then there is a subsequence such that $d_{i_2}/d_{i_1},d_{i_3}/d_{i_2},\dotsc$ tends to infinity, implying 
$\bar{d}(A)\leq m(\{d_{i_1},\dotsc,\})\leq m(\{1\})^k$
for any positive integer $k$. In fact our result is stronger:

\begin{theorem}\label{smpand} 
Suppose $d\geq 2$ and let $D_1,\dotsc,D_k\subset \R^+$ be  
arbitrary finite sets.
If the ratios $t_1/t_2,t_2/t_3,\dotsc,t_{k-1}/t_k$ tend to infinity, then
\begin{equation*}
m_{\R^d}(t_1\cdot D_1 \cup \dotsb \cup t_k\cdot D_k)\to \prod_{i=1}^k m_{\R^d}(D_i).
\end{equation*}
\end{theorem}

It is conceivable that there might be denser and denser $D$-avoiding
sets whose density approaches $m_{\R^d}(D)$ without
there being a $D$-avoiding set of density $m_{\R^d}(D)$. However, that is not
the case. We show that there is
a set which not just achieves this density, but whose measure cannot
be increased by an alteration on a bounded subset. Moreover, we show
that the constants $m_{\R^d}(D)$ can in principle be computed for any
finite set $D$. However, the high time complexity of our algorithm prohibits 
us from settling the question whether $m_{\R^2}(\{1\})<1/4$.

The principal tool of the paper is the so-called zooming-out lemma
stating that under the appropriate conditions we can ignore the
small-scale details of the measurable sets in question. In this
sense, it is similar to the celebrated Szemer\'edi regularity lemma.
The Szemer\'edi regularity lemma implies that for the purpose
of counting subgraphs every graph can be replaced
by a much smaller ``reduced graph'' \cite{komlos_simonovits_reglemm}.
The zooming-out lemma states that every measurable set can be
replaced by a ``zoomed-out set'' which captures some of information
about counting (by an appropriate integral) pairs of points that
are at a given distance away.

\section{The \texorpdfstring{$1$}{1}-dimensional case and the main idea}
Before delving into the proof of the results in $\R^d$ it is instructive
to examine the situation in $\Z$, for it is much simpler, of interest
on its own right, and illustrates some of the ideas used in the main results.

Throughout the paper we identify sets with their characteristic functions,
i.e., for a set $A$ we define $A(x)=1$ if $x\in A$ and $A(x)=0$ if
$x\not\in A$. In this section we use the notation $[a..b]$ to denote
the interval of the integers from $a$ to $b$, i.e., $[a..b]=\Z\cap[a,b]$. 

For a set $A\subset \Z$ define upper and lower densities by
\begin{equation*}
\bar{d}(A)=\limsup_{n\to\infty} \frac{\abs{A\cap [-n..n]}}{2n+1},\qquad 
\underbar{d}(A)=\liminf_{n\to\infty} \frac{\abs{A\cap [-n..n]}}{2n+1}.
\end{equation*}
The set $A$ is $D$-avoiding if $(A-A)\cap D=\emptyset$, where $A-A=\{a_1-a_2 : a_1,a_2\in A\}$ is the difference set of~$A$. Define $m(D)=\sup \bar{d}(A)$ 
where the supremum is over all $D$-avoiding sets. 

The simple-minded analogue of theorem~\ref{smpand} is false.
If $A$ is $\{1\}$-avoiding
set, then $A(x-1)+A(x)\leq 1$ and thus 
$2\abs{A\cap[-n..n]}\leq \sum_{k=-n}^{n+1} \bigl(A(x-1)+A(x)\bigr)\leq 2n+2$ showing
that $m(\{1\})\leq 1/2$. On the other hand, the set of even integers
shows that $m(\{1\})=1/2$. However, for every odd integer $t$ the
set of even integers shows that $m\bigl(\{1\}\cup t\cdot\{1\}\bigr)
=m(\{1,t\})=1/2$. 
This example also shows why the theorem~\ref{smpand} is itself false in 
$\R^1$. In $\R^1$ the integration of the inequality $A(x)+A(x+1)\leq 1$ 
yields $m(\{1\})\leq 1/2$. The set $\bigcup_{k\in\Z}[2k,2k+1)$ shows 
that $m(\{1\})=1/2$, and the same set shows that $m(\{1,t\})=1/2$
for every odd integer~$t$.

The version of theorem~\ref{smpand} that works in one dimension
involves excluding thickened sets, in order to avoid this kind of 
congruential obstacles. For a set $D\subset \Z$ we denote by
$D^k$ the $k$-neighborhood of $D$, i.e., $D^k=\{x\in\Z : \abs{x-y}\leq 
k\text{ for some } y\in D\}$.
\begin{theorem}\label{ideatheorem}
For every finite set $D_1\subset \Z$ there is a $k$ such that for
every finite non-empty set $D_2\subset \Z$ 
we have
\begin{equation*}
m\left(D_1\cup (t\cdot D_2)^k\right)<m(D_1)m(D_2)
\end{equation*}
for every positive integer $t$.
\end{theorem}
\begin{proof}
Denote $\diam D=\max_{d\in D} \abs{d}$.
Let $k$ be any even integer so that
$\diam D_1-k m(D_1)\leq -1$.

Suppose $A$ is $D_1\cup (t\cdot D_2)^k$-avoiding. Then the
set $A^{k/2}$ is $t\cdot D_2$-avoiding. To see that
suppose $x_1,x_2\in A^{k/2}$ is a pair of elements such that
$x_1-x_2\in t\cdot D_2$. By the definition of $A^{k/2}$ there are
$y_1,y_2\in A$ with $\abs{x_1-y_1}\leq k/2$ and $\abs{x_2-y_2}\leq k/2$.
By the triangle inequality $y_1-y_2\in (t\cdot D_2)^k$, which is a contradiction.

Write the set $A^{k/2}$ as a union of disjoint intervals
$A^{k/2}=[a_1..b_1]\cup[a_2..b_2]\cup\dotsb$ where
for no $i,j$ we have $b_i+1=a_j$. Each of these
intervals has length at least $k$. If $q$ is the smallest element of 
$D_2$, then none of these intervals has length exceeding $tq$, 
for $A^{k/2}$ is $\{tq\}$-avoiding. 
The density of $A^{k/2}$ does not exceed 
$m(t\cdot D_2)=m(D_2)$. The set $A$ is contained
$A^{k/2}$, so it suffices to bound the
density of $A$ on each of the intervals $[a_i..b_i]$.
By translating the interval $[a_i..b_i]$ it suffices
to consider the case $[a_i..b_i]=[0..n-1]$.

So, suppose $A'\subset[0..n-1]$ is $D_1$-avoiding and $\abs{A'}=s$.  Then
$\tilde{A}=A'+(n+\diam D_1)\Z=\bigcup_{t\in\Z} \bigl(A'+(n+\diam D_1)t\bigr)$ 
is $D_1$-avoiding because the copies of a $D_1$-avoiding set $A'$ 
are too far from
each other for there to be elements $x,y$ in different
copies such that $x-y\in D_1$. 
Since $\tilde{A}$ has density $s/(n+\diam D_1)\leq m(D_1)$, we
infer $s\leq m(D_1)n+\diam D_1$.

Now let us turn back to the proof of the theorem. 
For each interval $[a_i..b_i]$
the subintervals $[a_i..a_i+k/2-1]$ and $[b_i-k/2+1..b_i]$
do not meet~$A$.
Thus each interval in $A^{k/2}$ of length $n$ contains 
no more than  $m(D_1)(n-k)+\diam D_1\leq m(D_1)n-1\leq (m(D_1)-1/tq)n$ 
elements of $A$. 

Similarly no more than 
$m(t\cdot D_2)r+\diam(t\cdot D_2)=m(D_2)r+t\diam D_2$ 
elements belong to $A^{k/2}$ in any 
interval of length~$r$.
Let $n$ be an arbitrary positive integer. 
Consider $B=\bigcup_{[a_i..b_i]\subset[-n..n]} [a_i..b_i]$.
Since at most two intervals contain elements in $[-n..n]$,
but not contained in $[-n..n]$, we have
$\abs{B}\geq \abs{[-n..n]\cap A}-2tq$. 
Hence,
\begin{equation*}
\abs{[-n..n]\cap A}
\leq \abs{B}+2tq\leq
\bigl(m(D_2)(2n+1)+t\diam D_2\bigr)\bigl(m(D_1)-1/tq\bigr)+2tq.
\end{equation*}
Letting $n\to\infty$ we conclude
that $\bar{d}(A)\leq m(D_2)\bigl(m(D_1)-1/tq\bigr)<m(D_1)m(D_2)$.
\end{proof}

As remarked above the reason why theorem~\ref{smpand} fails
in the dimension one is because the largest $D$-avoiding set can be periodic
(in fact there is always a set of density $m(D)$ which is periodic as shown
by Cantor and Gordon \cite{cite:cantor_gordon_missingdistances}), and 
thus avoid many more distances than required of it. By the theorem
of Furstenberg, Katznelson and Weiss \cite{cite:fkw} this cannot
happen in higher dimensions because any periodic set has positive
density, and all sufficiently large distances occur in sets of
positive density. So, it is not surprising that in the higher dimensions it 
becomes possible to carry out a proof very similar in spirit to the 
proof of theorem~\ref{ideatheorem} above, but technically more complicated.

The approach employed in this paper is rooted in the proof of 
Bourgain \cite{cite:bourgain} of the Furstenberg-Katznelson-Weiss theorem.

\section{Notation}\label{sec_notation}
Throughout the rest of the paper the dimension $d\geq 2$ is going to 
be fixed, so we will often omit the dependency on $d$ from our notation.

For a measurable set $A\subset \R^d$ the notation $\meas{A}$ denotes
the measure of $A$. The notation $Q(x,r)$ denotes the open axis-parallel cube of
side length~$r$ centered at the point $x$. 

For a set $A\subset \R^d$ and a bounded domain $\Omega\subset \R^d$
the density of $A$ on $\Omega$ is
\begin{equation*}
d_{\Omega}(A)=\frac{\meas{A \cap \Omega}}{\meas{\Omega}}.
\end{equation*}
The upper and lower limit densities of $A$ are
\begin{equation*}
     \bar{d}(A)=\limsup_{R\to\infty} d_{Q(0,R)}(A),\qquad
\underbar{d}(A)=\liminf_{R\to\infty} d_{Q(0,R)}(A).
\end{equation*}
Whenever $\bar{d}(A)=\underbar{d}(A)$ we write 
$d(A)=\bar{d}(A)=\underbar{d}(A)$. 
Note that we measure the densities with respect to cubes, and not balls as
it is usually done. 
  Whereas, in general these densities might be different, corollary 
\ref{density_indep} below implies that
our results do not depend on the kind of density chosen, and the proofs
are cleaner for the density measured on cubes since
there are fewer edge effects one needs to worry about. 
The advantage of using cubes centered at the origin lies in
less cluttered notation. However, since the properties we
consider in this paper are translation-invariant, we incur
no loss of generality.

Being interested in the largest $D$-avoiding sets, we define
\begin{equation*}
m(D)=\sup_{A\text{ is $D$-avoiding}} \bar{d}(A).
\end{equation*}
More generally, we will be looking at the properties of sets that are more general 
than the property of being $D$-avoiding.  So, we let $\Msr(\R^d)$ denote the family
of all the measurable subsets of $\R^d$ and call a function $P\colon \Msr(\R^d) \to \{0,1\}$
a \emph{property}. If $P(A)=1$, we say that $A$ has the property $P$,
and if $P(A)=0$, we say that $A$ does not have it. We define 
\begin{equation*}
m(P)=\sup_{P(A)=1}\bar{d}(A),\qquad
m_{\Omega}(P)=\sup_{P(A)=1} d_{\Omega}(A).
\end{equation*}
For a property $P$ and a real number $t>0$ the property $t\cdot P$ is the
property that holds for $A$ precisely when the property $P$ holds for $(1/t)\cdot A$.
This is in agreement with the definition of $t\cdot D$-avoiding set as a set $A$
such that $(1/t)\cdot A$ is $D$-avoiding. Note that the function $m$ is 
scale-invariant: for every $t>0$ we have $m(t\cdot P)=m(P)$.

If $P_1$ and $P_2$ are two properties, then $P_1 \AND P_2$ denotes the
property asserting that both $P_1$ and $P_2$ hold, i.e., $(P_1 \AND P_2)(A)=
P_1(A)P_2(A)$. In particular, if $P_1$ and $P_2$ are the properties
of being $D_1$- and $D_2$-avoiding respectively, then $P_1 \AND P_2$
is the property of being $D_1\cup D_2$-avoiding. 

\section{Supersaturable properties}
In this section we prove basic theorems about a class of properties for 
which the analogue of theorem~\ref{smpand} holds. 

As explained in the introduction, the crucial tool is the ability to ignore the fine 
details of the sets. 
The intuition here is that given a set $A$ and
a large real number $t$ in order to understand whether the set has points 
which are at distance $t$ apart we should zoom-out away from the set $A$ and
 look at a
scale comparable to $t$. If we think of the set $A$ as colored black
on the otherwise white background, then the very fine details of $A$ will 
blur into some shade of gray. The zooming-out lemma says that for our purposes 
if the shade is not too light, then we can treat gray points as if they were black.

More formally, for each $\delta>0$ and $\epsilon>0$ we define a zooming-out 
operator $\Zm_\delta(\epsilon)$ acting on $\Msr(\R^d)$ by
\begin{equation*}
\Zm_\delta(\epsilon) A = \left\{x \in \R^d :
\abs{Q(x,\delta)\cap A}>\epsilon\abs{Q(x,\delta)}\right\}.
\end{equation*}
One can think of the zooming-out operator as the replacement for the operation 
of thickening sets $A\mapsto A^{k/2}$ in the integers. In the sequel we
use the following easy properties of the zooming-out operator which we
now state.
\renewcommand*{\theenumi}{\alph{enumi}}
\begin{lemma}\label{Zmproplemma}
\begin{enumerate}
\item\label{Zmpropone} $\Zm_{\delta}(\epsilon)A+Q(0,(t-1)\delta)
\subset \Zm_{t\delta}(t^{-d}\epsilon)A$ for
any $t\geq 1$.
\item\label{Zmproptwo} $\Zm_{\delta_1}(\epsilon_1)\Zm_{\delta_2}(\epsilon_2)A
\subset
\Zm_{\delta_1+\delta_2}\left(\epsilon_1\epsilon_2
\frac{\delta_1^d \delta_2^d}{(\delta_1+\delta_2)^d\min(\delta_1,\delta_2)^d}\right)A$.
\end{enumerate}
\end{lemma}
\begin{proof}
The claim \ref{Zmpropone} is clear, so we show \ref{Zmproptwo}. If
$x\in \Zm_{\delta_1}(\epsilon_1)\Zm_{\delta_2}(\epsilon_2)A$, then
\begin{equation*}
\epsilon_1\delta_1^d\leq \int \Zm_{\delta_2}(\epsilon_2)A(y)Q(x,\delta_1)(y)\,dy.
\end{equation*}
Since for all $y$ we have
\begin{equation*}
\epsilon_2\delta_2^d \Zm_{\delta_2}(\epsilon_2)A(y)
\leq \int A(z)Q(y,\delta_2)(z)\,dz,
\end{equation*}
it follows that
\begin{equation*}
\epsilon_1\epsilon_2\delta_1^d\delta_2^d\leq \int A(z) \meas{Q(z,\delta_2)\cap
Q(x,\delta_1)}\,dz\leq \min(\delta_1,\delta_2)^d
\meas{A\cap Q(x,\delta_1+\delta_2)}.
\end{equation*}
Therefore $x\in \Zm_{\delta_1+\delta_2}\left(\epsilon_1\epsilon_2
\frac{\delta_1^d \delta_2^d}{(\delta_1+\delta_2)^d\min(\delta_1,\delta_2)^d}\right)A$ as desired.
\end{proof}

We say that a property $P$ is \emph{supersaturable} if there is a function
$I_P\colon \Msr(\R^d)\to [0,\infty]$ such that the following seven
conditions are satisfied:
\renewcommand*{\theenumi}{\Roman{enumi}}
\begin{enumerate}
\item \label{sups_nontriviality} $0<m(P)$.
\item \label{sups_monotone} $I_P(A)$ is monotone nondecreasing and $P$ is monotone, i.e., $I_P(A)\geq I_P(B)$ and $P(A)\leq P(B)$ if $A\supset B$.
\item \label{sups_positivity} $I_P(A)>0$ implies that $A$ does not have 
the property $P$.
\item \label{sups_translation} Both $P$ and $I_P$ are translation-invariant: $P(A)=P(A+x)$ and $I_P(A)=I_P(A+x)$ for every $x\in \R^d$.
\item \label{sups_locality} There is a real number, which we denote by $\diam(P)$, such that
if $A_1$ and $A_2$ are sets which are at distance at least $\diam(P)$ away
from each other, then $I_P(A_1\cup A_2)\geq I_P(A_1)+I_P(A_2)$ and $A_1\cup A_2$ has
the property $P$ iff both $A_1$ and $A_2$ have the property $P$ .
\item \label{sups_zoomability} There is an $\epsilon>0$ and a strictly positive
function $f$, such that if $\Zm_\delta(1-\epsilon)A$ does not
have the property $P$, then $I_P(A)\geq f(\delta)$.
\item (Zooming-out lemma) \label{sups_zoomingout} If $A\subset Q(0,R)$ then 
$I_P(A)\geq g_P(\epsilon) I_P(\Zm_\delta(\epsilon)A)-
h_P(\epsilon,\delta)R^d$, where $g_P$ is positive and 
$h_P(\epsilon,\delta)\to 0$ as $\delta\to 0$ for any fixed $\epsilon>0$.
\setcounter{saveenum}{\value{enumi}}
\end{enumerate}
We call $I_P$ a \emph{saturation function} for the property $P$. 
An example of supersaturable property to keep in mind is the property of being $\{1\}$-avoiding, for which
the saturation function can be chosen to be 
$I(A)=\iint A(x)A(x+y)\,d\sigma(y)\,dx$ where
$\sigma$ is the uniform measure on the unit circle, and here for the second 
time we use the convention that a set $A$ is identified with the 
characteristic function of~$A$. In this example, with the exception of the zooming-out lemma all the conditions are not hard to check, and the zooming-out lemma will be proved in section~\ref{zoomout_sect}.
More generally in theorem~\ref{distsatur} we will show
that the property of being $D$-avoiding is an example of a supersaturable property.
The proof of theorem is independent of the results in this section, and
might be read before this section.

The motivation for the definition of the supersaturable properties is that not 
only $I_P(A)>0$ implies that $A$ does not have the property
$P$, but also $\bar{d}(A)>m(P)$ implies that $I_P(A)>0$. The latter
statement is the content of the following lemma.
\begin{saturationlemma}
Let $P$ be a supersaturable property. 
For every $\epsilon_1,\epsilon_2>0$ there
is a constant $c=c(\epsilon_1,\epsilon_2)>0$
such that for any $R>0$ there is
$\delta_0=\delta_0(\epsilon_1,\epsilon_2,R)$
such that the following holds.
For any $\delta\leq\delta_0$ and any measurable set 
$A\subset \R^d$ if
\begin{equation*}
d_{Q(0,R)}(\Zm_\delta(\epsilon_1)A)>m_{Q(0,R)}(P)(1+\epsilon_2),
\end{equation*}
then
\begin{equation*}
I_P(A)\geq c R^d.
\end{equation*}
In particular, $A$ does not have the property $P$.

Moreover, $\delta_0(\epsilon_1,\epsilon_2,R)$ is a monotone
non-decreasing function of $R$ for any fixed $\epsilon_1,\epsilon_2$.
\end{saturationlemma}
Before proving the supersaturation lemma, we need two lemmas. The first
lemma shows that the
rate of convergence in the definition of $m(P)$ cannot be too slow,
whereas the second lemma assures us that we need not to worry about
small values of $R$.
{\piccaption{Tiling $T$}\parpic[r]{\includegraphics[scale=0.95]{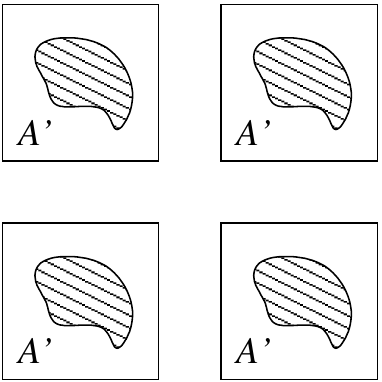}}\begin{lemma}\label{ratefact}
Let $P$ be a property satisfying the conditions \ref{sups_monotone}, \ref{sups_translation} and \ref{sups_locality}. If $A$ has the property $P$, then
\begin{equation*}
m(P)\geq d_{Q(0,R)}(A)\Big/\left(1+\frac{\diam P}{R}\right)^d.
\end{equation*}
\end{lemma}
\begin{proof}
Set $A'=A\cap Q(0,R)$. Then the tiling $T=A'+(R+\diam P)\Z^d$ has the property $P$
because the distance between the translates of
$A'$ is $\diam P$ and $A'$ has the property $P$. Since $m(P)\geq d(T)$,
the lemma follows.
\end{proof}}
\begin{lemma}\label{Rtozero}
Let $P$ be a property satisfying conditions \ref{sups_nontriviality},
\ref{sups_monotone} and
\ref{sups_translation}. Then $\lim_{r\to 0} m_{Q(0,r)}(P)=1$.
\end{lemma}
\begin{proof}
Assume the contrary. We will show there is no set of positive measure 
with property $P$, contradicting condition \ref{sups_nontriviality}.
Suppose there is a set $A$ of positive measure with property $P$. 
By the Lebesgue density theorem there is a point $p$ such that
$d_{Q(p,r)}(A)$ tends to $1$ as $r$ tends to $0$. By condition 
\ref{sups_translation} we may assume that $p=0$. 
Then the set $Q(0,r)\cap A$ is a subset of $Q(0,r)$ having property
$P$. Since the density of this set tends to $1$ as $r$ tends to zero 
we have reached a contradiction. 
\end{proof}
\begin{proof}[Proof of supersaturation lemma]
Since the condition of the lemma refers only to the set
$Q(0,R)\cap\Zm_\delta(\epsilon_1)A$ we can assume without any loss
of generality that $A\subset Q(0,R+2\delta)$. By lemma~\ref{Rtozero}
for every $\epsilon_2$ there is $R_{\text{min}}(\epsilon_2)>0$ 
such that if $R\leq R_{\text{min}}(\epsilon_2)$,
then $m_{Q(0,R)}(P)>1/(1+\epsilon_2)$. Thus
if $R\leq R_{\text{min}}(\epsilon_2)$, then
the premise of the supersaturation lemma
cannot hold since no set can have density 
$m_{Q(0,R)}(P)(1+\epsilon_2)>1$. Hence we can assume
that $R\geq R_{\text{min}}(\epsilon_2)>0$ throughout
the proof.

In the course of the proof of the supersaturation lemma we will prove 
following three statements:
\begin{itemize}
\item $\Lem(\epsilon_1,\epsilon_2)$ is the statement that the supersaturation
lemma holds for some specific $\epsilon_1$ and~$\epsilon_2$.
\item $\Lem'(\epsilon_2)$ is the statement that if
$A\subset Q(0,R)$ with $d_{Q(0,R)}(A)\geq m_{Q(0,R)}(P)(1+\epsilon_2)$, then the inequality
$I_P(A)\geq c' R^d$ holds with $c'=c'(\epsilon_2)>0$.
\item $\WkLem(1-\epsilon_T,\epsilon_2)$ is the statement
that for $\epsilon_1=1-\epsilon_T$
the conditions of the supersaturation lemma imply
the weaker conclusion in which the constant
$c$ is allowed to depend not only on $\epsilon_2$ but also on $\delta$.
Here $\epsilon_T$ is a positive number which depends only
on the property $P$.
\end{itemize}
First, we will establish $\WkLem(1-\epsilon_T,\epsilon_2)$
for every $\epsilon_2>0$. Then we will show that $\Lem'(\epsilon_2)$
implies $\Lem(\epsilon_1,\epsilon_2)$ for any $\epsilon_1>0$.
Finally, we will demonstrate that $\Lem(\epsilon_2 \epsilon_T m_{Q(0,R)}(P)/4,(1+\epsilon_T/8)\epsilon_2)$ and 
$\WkLem(1-\epsilon_T,\epsilon_2/2)$ together imply $\Lem'(\epsilon_2)$. 
Since for $\epsilon_2\geq 1/m_{Q(0,R)}(P)$ the $\Lem'(\epsilon_2)$
is vacuously true, all of
these imply $\Lem'(\epsilon_2)$ for all $\epsilon_2>0$ by induction
on $\bigl\lceil \log_{1+\epsilon_T/8} \frac{1}{\epsilon_2} \bigr\rceil$.
Then the proof will be complete.

\proofpart{$\WkLembf(1-\epsilon_T,\epsilon_2)$}We let 
$\epsilon_T$
to be the $\epsilon$ whose existence is postulated in the condition \ref{sups_zoomability}. We set $\delta_0=\diam P$.
Choose $R'$ so large that 
\begin{equation*}
\left(\frac{R'}{R'-3\diam P}\right)^d\leq 
\min\left(1+\frac{m_{Q(0,R)}(P)\epsilon_2}{4},
\frac{1+\epsilon_2/2}{1+\epsilon_2/3}
\right).
\end{equation*}
If $R\leq \frac{8d}{m_{Q(0,R)}(P)\epsilon_2}R'$ then since
$\Zm_\delta(1-\epsilon_T)(A)$ does not have
the property $P$, the condition 
\ref{sups_zoomability} tells us $I_P(A)\geq f(\delta)\geq 
c(\delta,\epsilon_2)R^d$ provided
$c$ is chosen small enough. So, assume
$R>\frac{8d}{m_{Q(0,R)}(P)\epsilon_2}R'$. Let $k=
\floor{R/R'}$. Let $\mathcal{C}_1$
be a collection of $k^d$ disjoint
cubes inside $Q(0,R)$ of side length $R'$ each.
Let $\mathcal{C}_2=
\{Q(x,R'-3\diam P) : Q(x,R')\in \mathcal{C}_1\}$.
Then 
\begin{align*}
d_{\bigcup \mathcal{C}_2}(\Zm_\delta(1-\epsilon_T)A)&\geq
\frac{
\meas{Q(0,R)}d_{Q(0,R)}(\Zm_\delta(1-\epsilon_T)A)
-\meas{Q(0,R)\setminus \bigcup \mathcal{C}_2}}{\meas{\bigcup \mathcal{C}_2}}\\
&=1-
\bigl(1-d_{Q(0,R)}(\Zm_\delta(1-\epsilon_T)A)\bigr)
\frac{R^d}{\meas{\bigcup \mathcal{C}_1}}\cdot
\frac{\meas{\bigcup \mathcal{C}_1}}{\meas{\bigcup \mathcal{C}_2}}\\
&\geq
1-
\bigl(1-d_{Q(0,R)}(\Zm_\delta(1-\epsilon_T)A)\bigr)
\frac{R^d}{(R-R')^d}\cdot \left(\frac{R'}{R'-3\diam P}\right)^d\\
&\geq 1-
\bigl(1-d_{Q(0,R)}(\Zm_\delta(1-\epsilon_T)A)\bigr)
\frac{1}{1-\frac{m_{Q(0,R)}(P)\epsilon_2}{8}}\cdot 
\left(1+\frac{m_{Q(0,R)}(P)\epsilon_2}{4}\right)\\
&\geq 
1-
\bigl(1-d_{Q(0,R)}(\Zm_\delta(1-\epsilon_T)A)\bigr)
\left(1+\frac{m_{Q(0,R)}(P)\epsilon_2}{2}\right)
\end{align*}
where in the last inequality we used that
$m_{Q(0,R)}(P)\epsilon_2<1$. Thus
\begin{align*}
d_{\bigcup \mathcal{C}_2}(\Zm_\delta(1-\epsilon_T)A)&\geq
d_{Q(0,R)}(\Zm_\delta(1-\epsilon_T)A)-\frac{m_{Q(0,R)}(P)\epsilon_2}{2}\\
&\geq m_{Q(0,R)}(P)(1+\epsilon_2/2)\\&\geq m(P)(1+\epsilon_2/2).\\
\intertext{From lemma~\ref{ratefact}, and the choice of $R'$
we get}
d_{\bigcup \mathcal{C}_2}(\Zm_\delta(1-\epsilon_T)A)&\geq m_{Q(0,R'-3\diam P)}(P)(1+\epsilon_2/3)
\end{align*}
Let $\mathcal{C}_3=\{Q(x,R'-3\diam P)\in\mathcal{C}_2 : 
d_{Q(x,R'-3\diam P)}(\Zm_\delta(1-\epsilon_T)A)>m_{Q(0,R'-3\diam P)}(P)(1+\epsilon_2/6)\}$. Set $n=\abs{\mathcal{C}_3}$. Then
\begin{equation}\label{presumone}
k^d m_{Q(0,R'-3\diam P)}(P)(1+\epsilon_2/3)\leq n+(k^d-n)m_{Q(0,R'-3\diam P)}(P)(1+\epsilon_2/6).
\end{equation}
Since by lemma~\ref{ratefact}
\begin{align*}
m_{Q(0,R'-3\diam P)}(P)(1+\epsilon_2/6)
&\leq m(P)\left(1+\frac{3\diam P}{R'-3\diam P}\right)^d
(1+\epsilon_2/6)\\
&\leq m_{Q(0,R)}(P)(1+\epsilon_2/2)<1,
\end{align*}
the inequality~\eqref{presumone} implies
\begin{equation*}
n\geq k^d\frac{m_{Q(0,R'-3\diam P)}(P)\epsilon_2/6}{1-m_{Q(0,R'-3\diam P)}(P)(1+\epsilon_2/6)}.
\end{equation*}
Since $\delta\leq \delta_0=\diam P$
we have
$Q(x,R'-3\diam P)\cap 
\Zm_\delta(1-\epsilon_T)A\subset 
\Zm_\delta(1-\epsilon_T)\bigl(Q(x,R'-\diam P)\cap A\bigr)$.
Therefore if $Q(x,R'-3\diam P)\in \mathcal{C}_3$,
then the condition~\ref{sups_zoomability} implies
$I_P(Q(x,R'-\diam P)\cap A)\geq f(\delta)$.
Since $Q(x_1,R'-\diam P)$ and $Q(x_2,R'-\diam P)$
are at distance at least $\diam P$ for distinct
$Q(x_1,R'-3\diam P),Q(x_2,R'-3\diam P)\in \mathcal{C}_3$,
we can apply the condition~\ref{sups_locality}
to deduce
\begin{equation*}
I_P\left(A\cap \bigcup_{Q(x,R'-3\diam P)\in \mathcal{C}_3} Q(x,R'-\diam P)\right)\geq c(\delta,\epsilon_2) R^d.
\end{equation*}
The monotonicity condition \ref{sups_monotone}
allows us to conclude that $I_P(A)\geq c(\delta,\epsilon_2) R^d$.

\proofpart{$\Lembf'(\epsilon_2)$ implies
$\Lembf(\epsilon_1,\epsilon_2)$} Suppose a set $A$ satisfies
conditions of $\Lem(\epsilon_1,\epsilon_2)$. Then the zooming-out lemma
and $\Lem'(\epsilon_2)$ tell us that
\begin{equation*}
I_P(A)
\geq g(\epsilon_1) I_P(\Zm_\delta(\epsilon_1) A)-c_6 h(\epsilon_1,\delta) 
(R+\delta)^d
\geq g(\epsilon_1) c'(\epsilon_2) R^d-c_6 h(\epsilon_1,\delta) (R+\delta)^d.
\end{equation*}
If $\delta$ small enough, we obtain that $I_P(A)\geq c(\epsilon_1,\epsilon_2)R^d$.

\proofpart{$\Lembf(\epsilon_2 \epsilon_T m_{Q(0,R)}(P)/4,(1+\epsilon_T/8)\epsilon_2)$ and 
$\WkLembf(1-\epsilon_T,\epsilon_2/2)$ imply $\Lembf'(\epsilon_2)$}%
With hindsight we set $\epsilon_1=\epsilon_2 \epsilon_T m_{Q(0,R)}(P)/4$.
Condition~\ref{sups_nontriviality}
asserts that $m(P)>0$ ensuring that $\epsilon_1>0$.
Recall that $R\geq R_{\text{min}}=R_{\text{min}}(\epsilon_2)$ and
let 
\begin{equation*}
\delta=\min\bigl(\diam P,R_{\text{min}}\epsilon_1/25d,
\delta_0(\epsilon_1,(1+\epsilon_T/8)\epsilon_2,R_{\text{min}})\bigr).
\end{equation*}
Suppose we have a set $A$ satisfying the conditions of $\Lem'(\epsilon_2)$.
If $A$ also satisfies the conditions of
$\WkLem(1-\epsilon_T,\epsilon_2/2)$, then $I_P(A)$ is as large as it should be, and we are done.
Hence, the conditions of $\Lem(1-\epsilon_T,\epsilon_2/2)$ do not hold.
Since $\delta\leq \diam P$, and $\delta_0$ in
$\WkLem(1-\epsilon_T,\epsilon_2/2)$ is equal to
$\diam P$, the only way in which the conditions of 
$\Lem(1-\epsilon_T,\epsilon_2/2)$ can fail is 
\begin{equation*}
d_{Q(0,R)}(\Zm_\delta(1-\epsilon_T)A)\leq m_{Q(0,R)}(P)(1+\epsilon_2/2).
\end{equation*}
Since the average density of $A$ is at least 
$m_{Q(0,R)}(P)(1+\epsilon_2)$ and the inequality above says that the density
of points that are centers of cubes of large density is no
more than $m_{Q(0,R)}(P)(1+\epsilon_2/2)$, there should be many
points that are centers of cubes with medium density $\epsilon_1$. 
For this we need to first relate $\meas{A}$ to $\meas{\Zm_\delta(\epsilon_1)A}$.
For that we need to allow for the edge effects due to averaging over the cube
of edge length $R+2\delta$ rather than $R$. Since 
$A\subset Q(0,R+2\delta)$,
\begin{align*}
\meas{A}&\leq \abs{A\cap Q(0,R-2\delta)}+
4d\delta (R+2\delta)^{d-1}\\&\leq
\left(\epsilon_1+(1-\epsilon_T)d_{Q(0,R)}(\Zm_\delta(\epsilon_1)A)+\epsilon_T d_{Q(0,R)}(\Zm_\delta(1-\epsilon_T)A)\right)R^d+5d\delta R^{d-1}.\\
\intertext{The definition of $A$ gives us}
\meas{A}&\geq m_{Q(0,R)}(P)(1+\epsilon_2)R^d.
\end{align*}
The two inequalities together yield
\begin{align*}
d_{Q(0,R)}(\Zm_\delta(\epsilon_1)A)&\geq \frac{m_{Q(0,R)}(P)\bigl(1-\epsilon_T+\epsilon_2(1-\epsilon_T/2)\bigr)-\epsilon_1-5d\delta/R}{1-\epsilon_T}.
\end{align*}
Our choice of $\epsilon_1$ and $\delta$, made in the beginning of the proof, 
assures
us that the left side is at least $m_{Q(0,R)}(P)(1+(1+\epsilon_T/8)\epsilon_2)$. Thus, we can apply
$\Lem\bigl(\epsilon_1,(1+\epsilon_T/8)\epsilon_2\bigr)$, 
and get the desired bound on
$I_P(A)$. This completes the proof of the final implication, and thus the
supersaturation lemma is proved.
\end{proof}
\picskip{0} One can combine the supersaturation
lemma with lemma~\ref{ratefact} to obtain a weak form
of supersaturation lemma which is easier to apply:
\begin{saturationlemmaalt}
Let $P$ be a supersaturable property. For every $\epsilon_1,\epsilon_2>0$
there are $\delta_0=\delta_0(\epsilon_1,\epsilon_2)>0$ and
$R_0=R_0(\epsilon_2)>0$ such that for any $\delta\leq \delta_0$ and
$R\geq R_0$ and any measurable set $A\subset \R^d$ if
\begin{equation*}
d_{Q(0,R)}(\Zm_\delta(\epsilon_1)A)>m(P)(1+\epsilon_2),
\end{equation*}
then
\begin{equation*}
I_P(A)\geq c R^d
\end{equation*}
for some constant $c=c(\epsilon_1,\epsilon_2)>0$
independent of $\delta$ and~$A$.
\end{saturationlemmaalt}
\begin{proof}
Choose $R_0$ to be large enough so that
$m_{Q(0,R)}(P)\leq m(P)(1+\epsilon_2/2)$ for $R\geq R_0$. Set
$\delta_0=\delta_0(\epsilon_1,\epsilon_2,R_0)$.
The monotonicity of $\delta_0(\epsilon_1,\epsilon_2,R)$ 
in the supersaturation lemma then insures that
any choice of $\delta\leq \delta_0$ and $R\geq R_0$
satisfies the conditions of the supersaturation lemma.   
\end{proof}

\begin{lemma}\label{lwrand}
If $P_1$ and $P_2$ are properties satisfying the conditions \ref{sups_translation} and \ref{sups_locality}, then $m(P_1 \AND P_2)\geq m(P_1)m(P_2)$.
\end{lemma}
\begin{proof}
Fix $\epsilon>0$. Take $R$ to be a large enough function of $\epsilon$.
Then pick a set $A_1$ with property $P_1$ such that $\bar{d}(A_1)\geq
m(P_1)(1-\epsilon)$. By averaging there is a cube $Q(x,R-\diam P_1)$ such that
$d_{Q(x,R-\diam P_1)}(A_1)\geq m(P_1)(1-2\epsilon)$. Then the proof of 
lemma~\ref{ratefact} shows existence of a periodic
set $A_1'$ with property $P_1$ of period $R$ with 
$d(A_1')\geq m(P_1)(1-3\epsilon)$.
Similarly, we can construct
a periodic set $A_2'$ with property $P_2$ with period $R$ 
and $d(A_2')\geq
m(P_2)(1-3\epsilon)$. Then averaging $d\bigl((A_1'+x)\cap A_2'\bigr)$ over
$x\in Q(0,R)$ yields existence of an $x_0$ such that
$d\bigl((A_1'+x_0)\cap A_2'\bigr)
\geq m(P_1)m(P_2)(1-3\epsilon)^2$. Since $\epsilon$
was arbitrary, the lemma follows.
\end{proof}
\begin{lemma}\label{ANDclosedness}
If $P_1$ and $P_2$ are any two supersaturable properties, then
so is $P_1 \AND P_2$. 
\end{lemma}
\begin{proof}
Let $I_{P_1}$ and $I_{P_2}$ be the saturation functions for $P_1$ and $P_2$.
Then $I_{P_1}+I_{P_2}$ is a saturation function for $P_1\AND P_2$.
The condition~\ref{sups_nontriviality}
follows from the lemma above. 
The conditions \ref{sups_monotone}, \ref{sups_positivity}, \ref{sups_translation} and \ref{sups_zoomingout} follow from the corresponding 
conditions for $P_1$ and $P_2$. For the conditions \ref{sups_locality} and
\ref{sups_zoomability} we can take $\diam (P_1 \AND P_2)=\max(\diam P_1,\diam P_2)$
and $\epsilon(P_1\AND P_2)=\min\bigl(\epsilon(P_1),\epsilon(P_2)\bigr)$ 
respectively.
\end{proof}
Now we are ready to derive a generalization of theorem~\ref{smpand}:
\begin{theorem}\label{ANDthm}
Suppose $P_1,\dotsc,P_n$ are supersaturable properties. Then
\begin{equation*}
m(t_1\cdot P_1 \AND \dotsb \AND t_n \cdot P_n)\to \prod_{i=1}^n m(P_i)
\end{equation*}
if for all $i\neq j$ the limit of $t_i/t_j$ is either $0$ or $\infty$.
\end{theorem}
\begin{proof}
The inequality $m(t_1\cdot P_1 \AND \dotsb \AND t_n \cdot P_n)\geq \prod_{i=1}^n m(P_i)$ 
follows from lemma~\ref{lwrand} and scale-invariance of $m$ by induction on~$n$.

For the proof of the opposite inequality we
permute $P_1,\dotsc,P_n$ and the corresponding
variables $t_1,\dotsc,t_n$ so that 
$t_{i+1}/t_i\to 0$ for all $i$. Furthermore, we scale
$t$'s so that $t_1=1$. Fix an arbitrary 
$\epsilon>0$. 
Let $\delta$ be the minimum of $\delta_0(\epsilon,\epsilon)$
over all the properties $P_1,\dotsc,P_{n-1}$,
where $\delta_0$ is as in the statement of the weak 
supersaturation lemma. Consider
any set $A$ with the property $t_1\cdot P_1 \AND \dotsb \AND t_n\cdot P_n$. 
Write $A_1=A$. The weak supersaturation lemma applied to this set and the property $P_1$ asserts that
\begin{equation*}
\bar{d}(\Zm_\delta(\epsilon)A_1) \leq m(P_1)(1+\epsilon).
\end{equation*}
For each point $x\in \Zm_\delta(\epsilon)A_1$ the set 
$Q(x,\delta)\cap A_1$ has the property $t_2\cdot P_2 \AND \dotsb \AND t_n \cdot P_n$.
Therefore, the set $A_2=(1/t_2)\cdot \bigl((A_1-x)\cap Q(0,\delta)\bigr)$ has the property 
$P_2 \AND (t_3/t_2)\cdot P_3 \dotsb \AND (t_n/t_2)\cdot P_n$. The set $A_2$
is contained in the cube $Q(0,\delta/t_2)$. Since $t_2\to 0$ we can assume that
$t_2$ is small enough so that we can apply the
weak supersaturation lemma to the set $A_2$ and property $P_2$ 
to get\vspace{-0.3em}\piccaption{Recursive zooming-in}\parpic[r]{\includegraphics[scale=1.03]{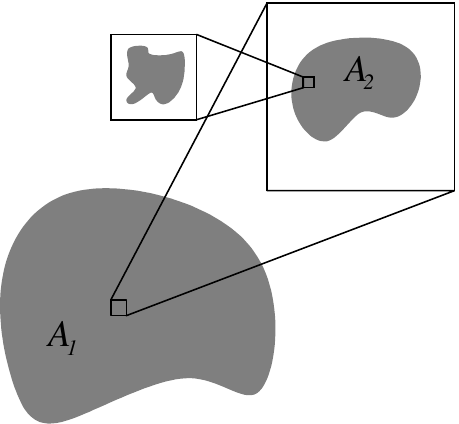}}%
\begin{equation*}
d_{Q(0,\delta/t_2)}(\Zm_\delta(\epsilon)A_2) \leq m(P_2)(1+\epsilon).
\end{equation*}
Repeating the argument, we eventually arrive at the inequalities
\begin{equation*}
d_{Q(0,\delta/t_{n-1})}(\Zm_\delta(\epsilon)A_{n-1}) \leq m(P_{n-1})(1+\epsilon)\\
\end{equation*}
and
\begin{align*}
d_{Q(0,\delta/t_n)} (A_n)&\leq m(P_n)(1+\epsilon).
\end{align*}
\picskip{3}These two inequalities mean that the density of $A_{n-1}$ 
on cubes of size $\delta$ is no more than $\epsilon$ except a set of
density no more than $m(P_{n-1})(1+\epsilon)$ on which the density is no more
than $m(P_n)(1+\epsilon)$. Hence, averaging implies that
\begin{align*}
d_{Q(0,\delta/t_{n-1})}(A_{n-1}) &\leq m(P_{n-1})m(P_n)(1+\epsilon)^2+\epsilon.
\end{align*}
Then by similarly unfolding the recursion, one arrives at the inequality
\begin{equation*}
\bar{d}(A_1)\leq \prod_{i=1}^n m(P_i) (1+\epsilon)^n+\mathcal{O}(\epsilon).
\end{equation*}
Since $\epsilon$ is arbitrary, this implies that $m(t_1\cdot P_1 \AND \dotsb \AND t_n\cdot P_n)\to\prod_{i=1}^n m(P_i)$.
\end{proof}

The definition of $m(P)$ leaves unclear whether there is ``a largest'' set with
property $P$ or there are larger and larger sets. If the property in question
is the property of not containing a copy of a finite subset in a given family,
then a largest set exists in a very strong sense.
\begin{definition}
A property $P$ is said to be \emph{finite} if there is a family $\mathcal{P}$
of finite sets such that $A$ has the property $P$ iff no set in $\mathcal{P}$ is a subset of~$A$. If in addition the 
diameter of sets in $\mathcal{P}$ is bounded,
then the property $P$ is said to be \emph{boundedly finite}.
\end{definition}
\begin{definition}
We call a measurable set $A\subset \Omega$ having property $P$ \emph{locally optimal} for the property $P$
with respect to a measurable set $\Omega$ if the following
condition holds for every bounded measurable set $S$: there is no measurable set $A'\subset\Omega$ with
property $P$ such that $A\cap (\R^d\setminus S)=A'\cap (\R^d\setminus S)$ 
such that $\meas{A'\cap S}>\meas{A\cap S}$.
If $\Omega=\R^d$, then we simply say that $A$ is locally optimal for~$P$.
\end{definition}
\begin{theorem}\label{locoptthm}
If $P$ is any boundedly finite 
supersaturable property and $\Omega$ is a measurable set, 
then there is a locally optimal set for~$P$ with respect to $\Omega$.
\end{theorem}
The proof of theorem~\ref{locoptthm}
requires an appropriate compactness result.
A characteristic function of any set lies in $L^\infty(\R^d)$,
which is a dual of $L^1(\R^d)$. The space $L^1(\R^d)$
induces a weak* topology
on $L^\infty(\R^d)$ which is the topology in which 
$f_1,f_2,\dotsc\to f$ when $\int f_k g\to\int fg$ as $k\to\infty$
for all $g\in L^1(\R^d)$.
\begin{lemma}\label{weakstarlim}
If $P$ is a finite supersaturable property, and
$A_1,A_2,\dotsc$ is a sequence of sets with property~$P$
whose characteristic functions converge
in the weak* topology of $L^\infty(\R^d)$.
Then there is a nonnegative
function $A\in L^\infty(\R^d)$ such that $A_1,A_2,\dotsc\to A$ in
the weak* topology, and 
$\supp A=\{x : A(x)>0\}$ has the property~$P$.
\end{lemma}
\begin{proof}
Since $A_1,A_2,\dotsc$ converge, they converge to some function, which
we will call~$A$.
The Lebesgue differentiation theorem states that
\begin{equation}\label{lebdiff}
\lim_{\delta\to 0}\frac{1}{\meas{Q(x,\delta)}}\int_{Q(x,\delta)}\abs{A(y)-A(x)}\,dy=0\qquad\text{for almost every $x$}.
\end{equation}
By setting $A$ to $0$ on a set of measure zero if necessary, we can
assume that this holds whenever $A(x)>0$ and $A$ is nonnegative. 
We will show that 
this modified function $A$ satisfies the conclusion of the lemma.
Suppose that on the contrary that
the set $\supp A$ lacked the
property~$P$. Then by finiteness of $P$
there would be a finite set $X=\{x_1,\dotsc,x_n\}\subset \supp A$
such that every set containing $X$ lacks~$P$.
Let $\epsilon=\min_{1\leq j\leq n} A(x_j)$.
Let $\epsilon_T$ be the $\epsilon$ whose
existence is postulated in the condition~\ref{sups_zoomability}.
By \eqref{lebdiff} there is $\delta$ such that for
every $1\leq j\leq n$ the
set
$\{y\in Q(x_j,\delta) : \abs{A(y)-A(x_j)}>\epsilon/4\}$
is of measure not exceeding $\tfrac{1}{4}\epsilon_T\epsilon\delta^d/24$.
Let $\delta'>0$ be any number small enough so that 
$(1-\delta'/\delta)^d>2/3$.
Let $Y_j=\{z\in Q(x_j,\delta-\delta') : \int_{Q(z,\delta')}\abs{A(y)-A(x_j)}\,dy
>\epsilon\delta'^d/4\}$. Since
\begin{align*}
\int_{Q(x_j,\delta)} \abs{A(y)-A(x_j)}\,dy&\geq \int_{Q(x,\delta-2\delta')} \frac{1}{\delta'^d}
\int_{Q(z,\delta')}\abs{A(y)-A(x_j)}\,dy\,dz\\
&\geq \meas{Y_j}\frac{\epsilon}{4},
\end{align*}
it follows that $\meas{Y_j}\leq \epsilon_T \delta^d/6$.
Let 
\begin{equation*}
W_j=\{y\in Q(x,\delta) : \abs{d_{Q(y,\delta')}(A_k)-\frac{1}{\delta'^d}\int_{Q(y,\delta')}A(z)\,dz}>\epsilon/4\}.
\end{equation*}
Choose $R$ to be so large that $Q(x_j,2\delta)\subset Q(0,R)$
for all $1\leq j\leq n$.
By the definition of weak* convergence for every $x$ we have
$\meas{A_{k}\cap Q(x,\delta')}=\int_{Q(x,\delta')} A_k(y)\,dy
\to \int_{Q(x,\delta')} A(y)\,dy$
as $k\to\infty$. So choose $k$ so large that 
$\abs{W_j}\leq \epsilon_T \delta^d/6$.
Thus for $y\in Q(x_j,\delta)\setminus(Y_j\cup W_j)$
we have $d_{Q(y,\delta')}(A_k)\geq \epsilon/2$.
Since $\meas{Y_j\cup W_j}\leq
\epsilon_T(\delta-\delta')^d/2$, we can also write this as 
$x_j\in\Zm_{\delta}(1-\epsilon_T/2)
\Zm_{\delta'}(\epsilon/2)\bigl(A_{k}\cap Q(0,R)\bigr)$. 
Let $t=\sqrt[d]{\frac{1-\epsilon_T/2}{1-\epsilon_T}}$.
By lemma~\ref{Zmproplemma}
\begin{equation*}
Q(x_j,(t-1)\delta)\subset \Zm_{t\delta}(1-\epsilon_T)
\Zm_{\delta'}(\epsilon/2)\bigl(A_{k}\cap Q(0,R)\bigr).
\end{equation*}
Since $x_j\in Q(x_j,(t-1)\delta)$, from the 
condition~\ref{sups_zoomability} we infer
$I_P\bigl(\Zm_{\delta'}(\epsilon/2)(A_{k}\cap Q(0,R)\bigr)>f(t\delta)$.
By the zooming-out lemma~\ref{sups_zoomingout}
we have
\begin{equation*}
I_P\bigl(A_{k}\cap Q(0,R)\bigr)\geq 
g_P(\epsilon/2) f(t\delta)-h_P(\epsilon/2,\delta')R^d.
\end{equation*}
Since $\delta'$ is independent of both $\epsilon$
and $\delta$, for sufficiently small $\delta'$ we 
would have that 
$I_P(A_{k})\geq I_P\bigl(A_{k}\cap Q(0,R)\bigr)>0$
contradicting the assumption that $A_{k}$ had
the property~$P$. The contradiction shows
that $\supp A$ has the property~$P$.
\end{proof}
\begin{proof}[Proof of theorem~\ref{locoptthm}]
Let $A_1,A_2\dotsc\subset \Omega$ be a sequence of sets,
each having the property $P$, such that
\begin{equation*}
d_{\Omega\cap Q(0,i)}(A_i)\geq m_{\Omega\cap Q(0,i)}(P)-2^{-i}.
\end{equation*}
We can and will assume that $A_i\subset \Omega\cap Q(0,i)$.
The Banach-Alaoglu theorem states that the closed ball in the dual of a 
Banach space is compact in weak* 
topology\cite[theorem~3.15]{cite:rudin_funanal}.
Thus there is a subsequence $A_{i_1},A_{i_2},\dotsc$ which converges
in weak* topology.
By lemma above there is a limit $A$ of the subsequence such that
the set $\supp A$ has the property~$P$. The set $\supp A$ is
the desired locally optimal set. Indeed, suppose that is not
so, and there are $R$, and $\epsilon>0$, and a set $A'\subset \Omega$ 
such that
$\meas{A'\cap Q(0,R)}\geq \meas{\supp A\cap Q(0,R)}+\epsilon$
and $A'\setminus Q(0,R)=\supp A\setminus Q(0,R)$. 
Since $P$ is boundedly finite, there is a $R'$ 
and a family of sets $\mathcal{P}$ of diameter
at most $R'$ each such that a set 
does not have the property $P$ precisely when the set
contains a member of $\mathcal{P}$.
Let $f$ be the characteristic function of 
$Q(0,R)\cap\supp A$.
By the definition of weak* convergence there are arbitrarily large $k$
so that 
\begin{align*}
&\left\lvert\int f(x)\bigl(A_{i_k}(x)-A(x)\bigr)\,dx\right\rvert\leq \epsilon/4,
\\\intertext{and}
&\left\lvert\int_{Q(0,R)} \bigl(A_{i_k}(x)-A(x)\bigr)\,dx\right\rvert\leq \epsilon/4.
\end{align*}
\piccaption{Set $\tilde{A}$}\parpic[r]{\begin{picture}(0,0)%
\includegraphics{optim.pdftex}%
\end{picture}%
\setlength{\unitlength}{3947sp}%
\begingroup\makeatletter\ifx\SetFigFont\undefined%
\gdef\SetFigFont#1#2#3#4#5{%
  \reset@font\fontsize{#1}{#2pt}%
  \fontfamily{#3}\fontseries{#4}\fontshape{#5}%
  \selectfont}%
\fi\endgroup%
\begin{picture}(1299,1299)(-11,-448)
\put(376,-136){\makebox(0,0)[lb]{\smash{{\SetFigFont{12}{14.4}{\rmdefault}{\mddefault}{\itdefault}{\color[rgb]{0,0,0}\small$A'\cap A_{i_k}$}%
}}}}
\put(526,-361){\makebox(0,0)[lb]{\smash{{\SetFigFont{12}{14.4}{\rmdefault}{\mddefault}{\itdefault}{\color[rgb]{0,0,0}\small$A_{i_k}$}%
}}}}
\put(564,164){\makebox(0,0)[lb]{\smash{{\SetFigFont{12}{14.4}{\rmdefault}{\mddefault}{\itdefault}{\color[rgb]{0,0,0}\small$A'$}%
}}}}
\end{picture}%
}
Let $\tilde{A}=\bigl(A'\cap Q(0,R)\bigr)\cup\Bigl(\supp A\cap A_{i_k}\cap
\bigl(Q(0,R+R')\setminus Q(0,R)\bigr)\Bigr)
\cup \bigl(A_{i_k}\setminus Q(0,R+R')\bigr)$.
Note that $\tilde{A}\cap Q(0,R+R')\subset A'\cap Q(0,R+R')$.

\picskip{5} If $\tilde{A}$ did not have the property~$P$,
then there would be a finite set $X\subset \tilde{A}$ such that
every set containing $X$ does not have the property~$P$.
By the definition of $R'$, we would have that either $X$
is a subset of either
$\bigl(A'\cap Q(0,R)\bigr)\cup\Bigl(\supp A\cap A_{i_k}\cap
\bigl(Q(0,R+R')\setminus Q(0,R)\bigr)\Bigr)$ or
$\Bigl(\supp A\cap A_{i_k}\cap
\bigl(Q(0,R+R')\setminus Q(0,R)\bigr)\Bigr)
\cup \bigl(A_{i_k}\setminus Q(0,R+R')\bigr)$. 
Since the former is a subset of $A'$ and the latter is 
a subset of $A_{i_k}$, we would reach a contradiction with
the assumption that $A'$ and $A_{i_k}$ both have the property~$P$.
Thus, $\tilde{A}$ has the property~$P$.

On the other hand, 
\begin{align*}
\meas{\tilde{A}}&=\meas{A'\cap Q(0,R)}+
\Bigl\lvert\supp A\cap A_{i_k}\cap
\bigl(Q(0,R+R')\setminus Q(0,R)\bigr)\Bigr\rvert
+\meas{A_{i_k}\setminus Q(0,R+R')}\\
&\geq \epsilon+\meas{\supp A\cap Q(0,R+R')\cap A_{i_k}}+
\meas{A_{i_k}\setminus Q(0,R+R')}\\
&= \epsilon+\int f(x) \bigl(A_{i_k}(x)-A(x)\bigr)\,dx+
\int_{Q(0,R+R')} \bigl(A(x)-A_{i_k}(x)\bigr)\,dx 
\\&\qquad+\meas{A_{i_k}\cap Q(0,R+R')}+\meas{A_{i_k}\setminus Q(0,R+R')}\\
&\geq \epsilon/2+\meas{A_{i_k}}.
\end{align*}
If $k$ was chosen large enough, we 
obtain $d_{\Omega\cap Q(0,i_k)}(\tilde{A})
\geq m_{\Omega\cap Q(0,i_k)}(P)-2^{-i_k}+i_k^{-d}\epsilon/2>
m_{\Omega\cap Q(0,i_k)}(P)$. The contradiction implies
that $\supp A$ is locally optimal.
\end{proof}
\begin{corollary}\label{density_indep}
If $P$ is any boundedly finite supersaturable property, then there is a set $A$ with
property $P$ such that for any open bounded set $\Omega$ 
\begin{equation*}
\lim_{t\to\infty} d_{t\cdot \Omega}(A)=m(P).
\end{equation*}
\end{corollary}
\begin{proof}
It follows from Whitney decomposition, for example, that we can 
write $\Omega$ as a union of countably many disjoint open cubes and a set of
measure zero, i.e., $\Omega=Z\cup \bigcup_{i\geq 0} Q(x_i,r_i)$ where $Z$ is 
of measure zero. Let $\epsilon>0$ be arbitrary and let $A$ be a locally optimal 
set for the property~$P$. Choose $n$ to be large enough so that $\meas{\bigcup_{i>n}
Q(x_i,r_i)}<\epsilon$.
 
By lemma~\ref{ratefact} the measure of $\bigl(t\cdot Q(x_i,r_i)\bigr)\cap A=Q(t x_i,t r_i)\cap A$ 
cannot exceed
$\meas{Q(t x_i,t r_i+\diam P)}m(P)$. By the local optimality of $A$ the measure of
$Q(t x_i,t r_i)\cap A$ cannot be any less than $\meas{Q(t x_i,t r_i-2\diam P)}m(P)$.
Hence
\begin{equation*}
\meas{Q(t x_i,t r_i)\cap A}=\meas{Q(tx_i,tr_i)}m(P)\bigl(1+\mathcal{O}(1/r_i t)\bigr).
\end{equation*}
Summing over $i$ with $i\leq n$ we obtain
\begin{equation*}
\bigl\lvert\meas{(t\cdot \Omega)\cap A}-\meas{t\cdot \Omega}m(P)\bigr\rvert<\epsilon+
\mathcal{O}\bigl(\tfrac{1}{t}\sum_{i\leq n} \tfrac{1}{r_i}\bigr).
\end{equation*}
Since $t$ goes to infinity and $\epsilon$ is arbitrary, the corollary follows.
\end{proof}

\section{Zooming-out lemma}\label{zoomout_sect}
In this section we establish that several properties
including the property of being $D$-avoiding are supersaturable.

We shall use Fourier transform on $\R^d$ which is defined via
\begin{equation*}
\hat{f}(\xi)=\int_{\R^d} f(x)e^{-2\pi i \langle x,\xi\rangle}\,dx,\qquad
\hat{\sigma}(\xi)=\int_{\R^d} e^{-2\pi i \langle x,\xi\rangle}\,d\sigma(x)
\end{equation*}
for a function $f$ or a Borel measure $\sigma$, respectively. For functions
$f,g\in L^1\cap L^\infty$ 
and a measure $\sigma\in \mathcal{M}(\R^d)$ the convolutions are defined
by $(f*g)(y)=\int f(y-x)g(x)\,dx$ and $(f*\sigma)=\int f(y-x)\,d\sigma(x)$,
which satisfy the following well-known identities
\begin{equation}\label{eq:idts}
\begin{aligned}
\widehat{f*g}(\xi)&=\hat{f}(\xi)\hat{g}(\xi),&\qquad&& 
\int f(x)g(x)\,dx&=\int \hat{f}(\xi)\hat{g}(-\xi)\,d\xi,\\
\widehat{f*\sigma}(\xi)&=\hat{f}(\xi)\hat{\sigma}(\xi),&&&
\int f(x)d\sigma(x)&=\int \hat{f}(-\xi)\hat{\sigma}(\xi)\,d\xi.
\end{aligned}
\end{equation}
\begin{definition}
A probability measure $\sigma\in\mathcal{M}(\R^d)$
with support $\supp \sigma$ 
is \emph{admissible} if $\sigma$ is symmetric
around $0$, has compact support, $0\not\in \supp \sigma$,
and $\hat{\sigma}(\xi)\to 0$ as $\abs{\xi}\to \infty$.
\end{definition}
We say that a set $A$ is  $\sigma$-avoiding
if there are no points $x,y\in A$ such that $x-y\in \supp \sigma$.
For example, the property of being $\{1\}$-avoiding in Euclidean
distance is the same as being $\sigma$-avoiding for $\sigma$ being
the surface measure on the unit sphere.  We can assume without loss
of generality that $\sigma$ is symmetric around $0$. Indeed if we let
$\sigma'(A)=\bigl(\sigma(A)+\sigma(-A)\bigr)/2$,
then being $\sigma'$-avoiding is same as being $\sigma$-avoiding.
Define the saturation function for the property of being $\sigma$-avoiding by 
$I_\sigma(A)=\int A(x)A(x+y)\,d\sigma(y)\,dx$. The saturation function
is well-defined by Tonelli's theorem.

Write $Q_\delta$ for the function $Q_\delta(x)=\delta^{-d} Q(0,\delta)(x)$.
\begin{lemma}\label{convlem} There is an absolute constant $c_1$ such that
if $\sigma\in \mathcal{M}(\R^n)$ is a probability measure, then for every $T>0$
\begin{equation*}
\left\lvert\iint f(x)g(x+y)\,d\sigma(y)\,dx-\iint f(x)(g*Q_\delta)(x+y)\,d\sigma(y)\,dx\right\rvert\leq \left(c_1\delta^{2}T^2+
\sup_{\abs{\xi}>T}2\abs{\hat{\sigma}(\xi)}
\right)\norm{f}_{L^2}\norm{g}_{L^2}.
\end{equation*}
\end{lemma}
\begin{proof}
Applying~\eqref{eq:idts} we obtain
\begin{align*}
\iint f(x) g(x+y)\,d\sigma(y)\,dx&=\int f(x) 
\int g(x+y)\,d\sigma(y)\,dx=\iint f(x)e^{-2\pi i\langle -x,\xi\rangle}\hat{g}(\xi)\hat{\sigma}(-\xi)\,dx\,d\xi\\&=\int \hat{f}(-\xi)\hat{g}(\xi)\hat{\sigma}(-\xi)\,d\xi
\end{align*}
Since $\abs{\widehat{Q_\delta}(\xi)-1}\leq c_1 \delta^2 \abs{\xi}^2$ and
$\abs{\widehat{Q_\delta}(\xi)}\leq \abs{\widehat{Q_\delta}(0)}=1$, we obtain
\begin{align*}
\left\lvert\iint f(x)\bigl(g(x+y)-(g*Q_\delta)(x+y)\bigr)\,d\sigma(y)\,dx\right\rvert
&=\left\lvert\int \hat{f}(-\xi)\hat{g}(\xi)\bigl(1-\widehat{Q_\delta}(\xi)\bigr)\hat{\sigma}(-\xi)
\,d\xi\right\rvert\\
&\leq c_1\delta^2 T^2\norm{\sigma} \int_{\abs{\xi}<T} \left\lvert\hat{f}(-\xi)\hat{g}(\xi)\right\rvert\,d\xi
\\ &\quad+\sup_{\abs{\xi}>T}2\abs{\hat{\sigma}(\xi)}\int_{\abs{\xi}\geq T} 
\left\lvert\hat{f}(-\xi)\hat{g}(\xi)\right\rvert\,d\xi
\end{align*}
Cauchy-Schwarz and Parseval imply that 
$\int \abs{\hat{f}(\xi)\hat{g}(\xi)}\leq \norm{\hat{f}}_{L^2}
\norm{\hat{g}}_{L^2}=\norm{f}_{L^2}\norm{g}_{L^2}$, completing
the proof.
\end{proof}
\begin{corollary}
If a probability measure $\sigma\in \mathcal{M}(\R^n)$ is admissible, then the property of being $\sigma$-avoiding satisfies the condition
\ref{sups_zoomingout}.
\end{corollary}
\begin{proof}
Suppose $A\subset Q(0,R)$. By the definition of $\Zm_\delta(\epsilon)A$ we
have
\begin{equation*}
\epsilon\Zm_\delta(\epsilon)A(x)\leq (A*Q_\delta)(x)
\end{equation*}
which implies
\begin{equation*}
\epsilon^2 I_\sigma\bigl(\Zm_\delta(\epsilon)A\bigr)\leq 
\iint (A*Q_\delta)(x)(A*Q_\delta)(x+y)\,d\sigma(y)\,dx.
\end{equation*}
Since $\sigma$ is symmetric around $0$, we have
\begin{equation*}
\iint (A*Q_\delta)(x)A(x+y)\,d\sigma(y)\,dx=
\iint A(x)(A*Q_\delta)(x+y)\,d\sigma(y)\,dx
\end{equation*}
and the lemma \ref{convlem} applied 
twice yields
\begin{align*}
&\left\lvert\iint (A*Q_\delta)(x)(A*Q_\delta)(x+y)\,d\sigma(y)\,dx-
\iint A(x)A(x+y)\,d\sigma(y)\,dx\right\rvert\\
&\qquad\leq \left\lvert\iint (A*Q_\delta)(x)(A*Q_\delta)(x+y)\,d\sigma(y)\,dx-
\iint (A*Q_\delta)(x)A(x+y)\,d\sigma(y)\,dx\right\rvert
\\&\qquad\quad+
\left\lvert\iint (A*Q_\delta)(x)A(x+y)\,d\sigma(y)\,dx-
\iint A(x)A(x+y)\,d\sigma(y)\,dx\right\rvert\\
&\qquad\leq \left(c_1\delta^{2}T^2+
\sup_{\abs{\xi}>T}2\abs{\hat{\sigma}(\xi)}
\right)\left(\norm{A}_{L^2}\norm{A*Q_\delta}_{L^2}+\norm{A}_{L^2}^2\right)
\end{align*}
Since $\norm{A}_{L^2}^2=\abs{A}\leq R^d$ and $
\norm{A*Q_d}_{L^2}=\norm{\hat{A} \widehat{Q_d}}_{L^2}
\leq \norm{\hat{A}}_{L^2}=\norm{A}_{L^2}$, it follows
that
\begin{equation*}
I_\sigma(A)\geq \epsilon^2 I_\sigma(\Zm_\delta(\epsilon)A)- 2\left(c_1\delta^{2}T^2+
\sup_{\abs{\xi}>T}2\abs{\hat{\sigma}(\xi)}
\right)R^d.
\end{equation*}
If we let $T=\delta^{-1/2}$, the condition
$\hat{\sigma}(\xi)\to 0$ as $\abs{\xi}\to\infty$
implies the condition \ref{sups_zoomingout}.
\end{proof}
With the zooming-out lemma in place we are ready to show supersaturability:
\begin{theorem}\label{distsatur}
If $\sigma\in \mathcal{M}(\R^n)$ is admissible, then the property of
being $\sigma$-avoiding is supersaturable.
\end{theorem}
\begin{proof}
The conditions \ref{sups_monotone}, \ref{sups_positivity}, \ref{sups_translation} are obvious.
The compact support of $\sigma$ implies the condition \ref{sups_locality}. 
Since $0\not\in \supp \sigma$ there is an $\epsilon>0$ such that
$Q(0,\epsilon)\cap \supp \sigma=\emptyset$. Then the set
$Q(0,\epsilon/2)+\diam(\supp \sigma) \Z^d$ has positive density and is
$\sigma$-avoiding. Thus the condition \ref{sups_nontriviality} is fulfilled.

Finally to verify the condition~\ref{sups_zoomability}
let $\epsilon=1/4$ and suppose $\Zm_\delta(1-\epsilon)A$ is not $\sigma$-avoiding.
Then there are $x_0,y_0\in \Zm_\delta(1-\epsilon)A$ such that
$x_0-y_0\in \supp \sigma$. Then for every $z\in Q(0,\delta/8d)$
the set $(A-x_0-z)\cap (A-y_0)\cap Q(0,\delta)$
has measure at least $\delta^d(1-2\epsilon-2d/8d)=\delta^d/4$.
Therefore
\begin{align*}
I_\sigma(A)&=\int A(x)A(x+y)d\sigma(y)\,dx\\
&\geq \iint_{x\in Q(x_0,\delta/8d)} A(x)A(x+y)\,dx\,d\sigma(y)\\
&\geq \frac{\delta^d}{4}\sigma\bigl(Q(y_0-x_0,\delta/8d)\bigr)
\end{align*}
which is positive since $y_0-x_0\in \supp \sigma$. 
\end{proof}
In particular since the surface measure on the unit sphere in $\R^d$ with 
$L^p$ norm for $1<p<\infty$ and $d\geq 2$ satisfies the condition
of the theorem, the property of being $\{1\}$-avoiding in $L^p$ is 
supersaturable property not only
for the usual Euclidean distance, but also in $L^p$ for 
$1<p<\infty$. By lemma \ref{ANDclosedness}
the property of being $D$-avoiding for a finite set $D\subset \R^+$ is also
of this form.

To avoid the false impression that the property
of being $\sigma$-avoiding is the only supersaturable
property, we demonstrate another class of natural
supersaturable properties. For symmetric
probability measures $\sigma_1,\sigma_2\in\mathcal{M}(\R^d)$ 
say  that a set $A$ is $\sigma_1 \OR \sigma_2$-avoiding 
if for every point $x\in A$ either there is no point $y\in
A$ such that $x-y\in\supp \sigma_1$ or
there is no point $y\in A$ such that 
$x-y\in \sigma_2$.
\begin{theorem}\label{ORsatur}
If $\sigma_1,\sigma_2$ are two admissible measures,
then the property of being $\sigma_1\OR \sigma_2$-avoiding is
supersaturable with the saturation function
$I_{\sigma_1 \OR \sigma_2}(A)=\iiint A(x)A(x+y_1)A(x+y_2)
\,d\sigma_1(y_1)\,d\sigma_2(y_2)\,dx$.
\end{theorem}
\begin{proof}
The conditions~\ref{sups_nontriviality} through \ref{sups_zoomability}
are checked in the same way as in the theorem~\ref{distsatur}.
We will show that \ref{sups_zoomingout} is satisfied.
Since $\sigma_1$ is a probability
measure, we have $\int A(x+y_1)\,d\sigma_1(y)\leq 1$ for every~$x$.
Therefore, lemma \ref{convlem} implies the inequality 
\begin{subequations}
\begin{align}\notag
&\left\lvert
\iiint A(x)A(x+y_1)\bigl(A(x+y_2)-(A*Q_{2\delta})(x+y_2)\bigr)
\,d\sigma_1(y_1)\,d\sigma_2(y_2)\,dx
\right\rvert
\\\label{eqORone} &\qquad\leq \left(4c_1\delta^{2}T^2+
\sup_{\abs{\xi}>T}2\abs{\hat{\sigma}_2(\xi)}
\right)\norm{A}_{L^2}^2.\\
\intertext{Similarly,}
&\notag\left\lvert
\iiint A(x)\bigl(A(x+y_1)-(A*Q_{2\delta})(x+y_1)\bigr)(A*Q_{2\delta})(x+y_2)
\,d\sigma_1(y_1)\,d\sigma_2(y_2)\,dx
\right\rvert\\\label{eqORtwo}&\qquad\leq
\left(4c_1\delta^{2}T^2+
\sup_{\abs{\xi}>T}2\abs{\hat{\sigma}_1(\xi)}
\right)\norm{A}_{L^2}^2.
\end{align}
\end{subequations}
Define translation operator by $(\mathcal{T}_x f)(z)=f(z-x)$.
Set $I'(f)=\iint f(y_1)f(y_2) \,d\sigma_1(y_1)\,d\sigma_2(y_2)$. Then
the inequalities \eqref{eqORone} and \eqref{eqORtwo} imply that
\begin{equation}\label{triagineq}
I_{\sigma_1\OR\sigma_2}(A)\geq 
\int A(x)I'(\mathcal{T}_x A*Q_{2\delta})\,dx-
\left(8c_1\delta^{2}T^2+
\sup_{\abs{\xi}>T}2\abs{\hat{\sigma}_1(\xi)}+
\sup_{\abs{\xi}>T}2\abs{\hat{\sigma}_2(\xi)}
\right)R^d.
\end{equation}
Since for every $y\in Q(x,\delta)$ we have 
$(\mathcal{T}_x A*Q_\delta)(z)\leq 
(2^d\mathcal{T}_y A*Q_{2\delta})(z)$ and
$I'$ is monotone, it follows that
\begin{equation}\label{transferineq}
\begin{aligned}
\int A(x)I'(\mathcal{T}_x A*Q_{2\delta})\,dx
&\geq 4^{-d}\delta^{-d} \iint_{y\in Q(x,\delta)}A(x)I'(\mathcal{T}_y A*Q_{\delta})\,dx\,dy\\
&=4^{-d}\delta^{-d} \iint_{y\in Q(0,\delta)}(\mathcal{T}_y A)(x)I'(\mathcal{T}_x A*Q_{\delta})\,dx\,dy\\
&=4^{-d}\int (A*Q_\delta)(x)I'(\mathcal{T}_x A*Q_{\delta})\,dx\\
&\geq 4^{-d}\epsilon^{3}I_{\sigma_1 \OR \sigma_2}(\Zm_{\delta}(\epsilon)A)
\end{aligned}
\end{equation}
If we set $T=\delta^{-1/2}$, the 
inequalities \eqref{triagineq} and \eqref{transferineq}
together imply the condition \ref{sups_zoomingout}.
\end{proof}

\section{Applications}
This section is devoted to two applications of the general results
proved above. 
\begin{theorem}\label{ORthm}
Let $\sigma_1,\sigma_2$ be a pair of admissible measures.
Let $P_1,P_2,P_1\OR t\cdot P_2$ denote the properties of being
$\sigma_1$-avoiding, $\sigma_2$-avoiding and $\sigma_1\OR 
t\cdot \sigma_2$-avoiding, respectively.
Then
\begin{equation*}
\lim_{t\to\infty} m(P_1\OR t\cdot P_2)=\max\bigl(m(P_1),m(P_2)\bigr).
\end{equation*}
\end{theorem}
\begin{theorem}\label{computthm}
There is an algorithm that given as input $\epsilon>0$ and a finite
set $D$ of distances outputs $m(D)$ with absolute error at most~$\epsilon$.
\end{theorem}

Before proving the theorem \ref{ORthm} we need some notation and a lemma.
If a set $A$ is $\sigma$-avoiding for admissible measure $\sigma$,
we set $F(A)=\{x\in \R^d : x-y\in \supp \sigma\text{ for some }y\in A\}$
and $S(A)=F(A)\cup A$.
Intuitively, if we try to enlarge $A$ to another
$\sigma$-avoiding set, then 
$F(A)$ is the sets which is forbidden by $A$ and
$S(A)$ is the set which is already ``occupied'' by~$A$.
Write $\diam \sigma=\max_{x\in\supp \sigma} \abs{x}$.
\begin{lemma}
Let $P$ be the property of being $\sigma$-avoiding.
For every $\sigma$-avoiding set $A$ we have $d_{Q(0,R)}(A)<m(P) 
\Bigl(d_{Q(0,R)}\bigl(S(A)\bigr)+(1+\diam \sigma/R)^d-1\Bigr)$.
\end{lemma}
\begin{proof}
We use the same trick that was used in the proof of lemma~\ref{ratefact}.
The set $A_1=\bigl(A\cap Q(0,R)\bigr)+(R+\diam \sigma)\Z^d$ is 
$\sigma$-avoiding and
$S(A_1)\subset \Bigl(\bigl(S(A)\cap Q(0,R)\bigr)\cup\bigl(Q(0,R+\diam \sigma)\setminus
Q(0,R)\bigr)\bigr)+(R+\diam \sigma)\Z^d$. Let $\alpha=d(A_1)/d(S(A_1))$. 
Since $d(S(A_1))\leq d_{Q(0,R)}(S(A))/(1+\diam \sigma/R)^d
+\bigl(1-1/(1+\diam \sigma/R)^d\bigr)$
and $d(A_1)\geq d_{Q(0,R)}(A)/(1+\diam \sigma/R)^d$,
it suffices to show that $\alpha\leq m(P)$.

Let $\gamma=1-d(S(A_1))$ be the proportion of $\R^d$ which
is not occupied yet. Choose a vector $x$ uniformly at random 
from $Q(0,R+\diam \sigma)$. For any set $X$
periodic with fundamental region $Q(0,R+\diam \sigma)$
we have that $E\bigl[d\bigl((X+x)\setminus S(A_1)\bigr)\bigr]=
\gamma d(X)$ where $E$ denotes the expectation. Let
$A_2=A_1\cup \bigl((A_1+x)\setminus S(A_1)\bigr)$.
Then $E\bigl[d(A_2)\bigr]=d(A_1)+\gamma d(A_1)$,
and $E\bigl[d\bigl(S(A_2)\bigr)\bigr]=d\bigl(S(A_1)\bigr)+\gamma d(S(A_1))$.
Hence $E\bigl[\alpha d\bigl(S(A_2)\bigr)-d(A_2)\bigr]=0$.
It follows that the set 
$\bigl\{x\in Q(0,R+\diam \sigma) : \alpha d\bigl(S(A_2)\bigr)-d(A_2)\leq 0\bigr\}$
has non-zero measure. In particular, it contains an element which is
not a period of the set~$A_1$. Thus, we can ensure 
$d(A_2)>d(A_1)$.

Similarly we can build an increasing sequence $A_1,A_2,A_3,\dotsc$ 
of $\sigma$-avoiding sets
such that $d(A_k)/d\bigl(S(A_k)\bigr)\geq \alpha$.
If the set $S(\bigcup_k A_k)$ had density~$1$, then we would be done, but that 
need not be the case. We use compactness lemma~\ref{weakstarlim} to circumvent this.

So, suppose $\alpha>m(P)$. Let $\mathcal{A}$ be
the family of all $\sigma$-avoiding sets $A\subset\R^d$ which are periodic
with the fundamental region $Q(0,R+\diam \sigma)$. Let $\mathcal{A}'$ be
those of them that satisfy $d(A)/d\bigl(S(A)\bigr)\geq \alpha$. 
Let $\beta=\sup_{A\in\mathcal{A}'} d(A)$.
Note that by the argument above the supremum is not achieved.
Let $A_1,A_2\dotsc\in\mathcal{A}'$ be a sequence such that $d(A_i)\to \beta$.

By passing to a subsequence if needed, assume that the sequence $A_1,A_2,\dotsc,$
converges in the weak* topology of $L^\infty(\R^d)$. 
By lemma~\ref{weakstarlim} there is a weak* limit $A$ of the sequence
such that $\supp A$ is $\sigma$-avoiding. Let $Y_i=
A_i\setminus \supp A$. We claim that
$d(Y_i)\to 0$ as $i\to\infty$. Suppose that there was a subsequence $Y_{i_1},Y_{i_2},\dotsc$
on which $d(Y_i)>\epsilon>0$. Banach-Alaoglu tells us that, by passing to a subsequence
again if needed, we can assume that $Y_{i_1},Y_{i_2},\dotsc$ converges
to some function~$Y$ in weak* topology. Since $\meas{\supp Y\cap \supp A}=0$,
we conclude that $\lim \int_{\supp Y} \bigl(A(x)-A_i(x)\bigr)\,dx<-\epsilon<0$ which contradicts
the definition of the weak* convergence. Thus, $d(Y_i)\to 0$.
Therefore, the sequence $A_i\cap \supp A$ converges to
$A$ in the weak* topology. 

Next we show that 
\begin{equation}\label{limorin}
d\bigl(F(\supp A)\setminus F(A_i)\bigr)\to 0\qquad\text{ as }i\to\infty.
\end{equation}
Pick an $\epsilon>0$. We will first
cover almost all of the set $F(\supp A)$ by cubes on which $F(\supp A)$
has density at least $1-\epsilon$. Then we will show that
$F(A_i)$ has density at least $1-3\epsilon$ on each of the cubes
provided $i$ is large.

Let $
\mathcal{Q}=\bigl\{Q(x,r) : 
d_{Q(x,r)}\bigl(F(\supp A)\bigr)>1-\epsilon\bigr\}$,
and let 
$\mathcal{C}$ be a family of all collections of
cubes from $\mathcal{Q}$ which are pairwise disjoint.
By Hausdorff maximum principle there is a maximal collection
$\mathcal{M}$ in $\mathcal{C}$. Then $W=F(\supp A)\setminus \bigcup \mathcal{M}$
is of measure null. Indeed, if $\meas{W}>0$ then by Lebesgue density
for almost every $x\in W$ we would have 
$\lim_{\delta\to 0} d_{Q(x,\delta)}(W)=1$, which implies
that there is $x$ and $\delta$ such that $d_{Q(x,\delta)}(W)>1-\epsilon$.
That contradiction shows that the desired covering exists.

Now let $Q(x_0,r)$ be any cube in the covering. Let $f(x)=
\int A(x+y)\,d\sigma(y)$. The function $f$ is defined almost
everywhere by Tonelli's theorem. Let 
$Z=\{x\in Q(x_0,r)\cap F(\supp A): f(x)=0\}$. The set
$Z$ is of measure null. Indeed, if $\meas{Z}>0$,
then by Lebesgue density theorem there would exist
an $x\in F(\supp A)$ such that $x\in \Zm_\delta(2/3)Z$ for all
sufficiently small $\delta$.
Let $y\in A$ be such that $\abs{x-y}\in\supp \sigma$.
Then since every point of $\supp A$ is a point of density,
there are arbitrarily small $\delta$ such that $y\in\Zm_\delta(2/3)(\supp A)$.
Thus, $\iint_{Q(x,\delta)} f(x)\,dx>0$. This contradicts the 
definition of $Z$ and so $\meas{Z}=0$.
Therefore,
there is a $\epsilon'$ such that the measure of $\{x\in
Q(x_0,r)\cap F(\supp A) : f(x)<\epsilon'\}$ does not
exceed $\epsilon r^d$. Therefore,
if $d_{Q(x_0,r)}(Y)\geq 3\epsilon$, then 
$\iint Y(x)A(x+y)\,d\sigma(y)\,dx>\epsilon \epsilon' r^d$.

Suppose there are arbitrarily large $i$'s
such that $F(A_i)$ has density less then $1-3\epsilon$ on
$Q(x_0,r)$. Let $Y_i=Q(x_0,r)\setminus F(A_i)$. Then
$\iint Y_i(x)A(x+y)\,d\sigma(y)\,dx>\epsilon \epsilon' r^d$.
Let $W_i=A_i\cap Q(x_0,r+2\diam P+1)$. Clearly,
$\iint Y_i(x)W_i(x+y)\,d\sigma(y)\,dx=0$. For small enough
$\delta$ lemma~\ref{convlem} implies that
\begin{equation}\label{convor}
\begin{aligned}
\iint Y_i(x)(W_i*Q_\delta)(x+y)\,d\sigma(y)\,dx&\leq \epsilon \epsilon' r^d/4\\
\iint Y_i(x)(A*Q_\delta)(x+y)\,d\sigma(y)\,dx&\geq \epsilon \epsilon' r^d/2.
\end{aligned}
\end{equation}
For any $\delta<1$ and for large enough $i$ we have
have $\abs{(W_i*Q_\delta)(x)-(A_i*Q_\delta)(x)}\leq
\frac{1}{9}\epsilon\epsilon'r^d (R+2\diam P+1)^{-d}$ for all $x\in Q(x,r+2\diam P)$ except a set of 
measure $\tfrac{1}{9}\epsilon'\epsilon'r^d$. Then
\begin{equation*}
\left\lvert\iint Y_i(x+y)(A_i*Q_\delta-W_i*Q_\delta)(x)\,d\sigma(y)\,dx\right\rvert\leq
\tfrac{2}{9}\epsilon'\epsilon'r^d,
\end{equation*}
which contradicts~\ref{convor}. Therefore, $F(W_i)$ has density
at least $1-3\epsilon$ on $Q(x_0,r)$ for all sufficiently large $i$.

By~\eqref{limorin} we get
\begin{equation*}
d\bigl(F(\supp A)\bigr)\leq \liminf_{i\to\infty} d\bigl(F(A_i)\bigr)
\leq \liminf_{i\to\infty} \left(\frac{1}{\alpha}-1\right)d(A_i)=\left(\frac{1}{\alpha}-1\right)\beta.
\end{equation*}
Since
\begin{equation*}
d(\supp A)\geq \frac{1}{\meas{Q(0,R+\diam \sigma)}}\int_{Q(0,R+\diam \sigma)}A(x)\,dx=\beta,
\end{equation*}
we conclude that
\begin{equation*}
\alpha d\bigl(F(\supp A)\bigr)\leq (1-\alpha)\beta+\alpha d(\supp A)\leq d(\supp A)
\end{equation*}
implying that $\supp A\in\mathcal{A}'$, which contradicts the assumption
that the supremum in the definition of $\beta$ is not achieved.
\end{proof}

\begin{proof}[Proof of theorem \ref{ORthm}]
The inequality $m(P_1\OR P_2)\geq \max\bigl(m(P_1),m(P_2)\bigr)$ is obvious.
Let $\epsilon>0$ be arbitrary. 
Let $R=\tfrac{5d}{\epsilon}\diam P_2$.
We will show that
$\limsup_{t\to 0} m_{Q(0,R)}\bigl((1/t)\cdot P_1\OR P_2\bigr)
\leq \max\bigl(m(P_1),m(P_2)\bigr)+\epsilon$. 

Suppose the contrary. Let $t_1,t_2\dotsc$ be a sequence
of $t$'s going to infinity for which $m_{Q(0,R)}\bigl((1/t_i)\cdot P_1\OR P_2\bigr)
\geq \max\bigl(m(P_1),m(P_2)\bigr)+\epsilon$.
Let $A_i$ be a locally optimal set for the property
$(1/t_i)\cdot P_1\OR P_2$. Let $A_i^1=
\{x\in A_i : \forall y\in A_i\ \abs{x-y}\not\in \tfrac{1}{t_i}\cdot
\supp \sigma_1\}$,
and $A_i^2=A_i\setminus A_i^1$. Note that
$A_i^1$ is $\sigma_1$-avoiding.
By passing to a subsequence
if necessary we can assume that the sequences $\{A_i^1\}_{i=1}^\infty$
and $\{A_i^2\}_{i=1}^\infty$ converge in weak*.
By lemma \ref{weakstarlim} there
is a limit $A^2$ of the sequence $\{A_i^2\}_{i=1}^\infty$ such that
$\supp A^2$ is $\sigma_2$-avoiding, and every point of $\supp A^2$
is a density point as in the Lebesgue density theorem.
Let $A^1$ be a limit of $\{A_i^1\}_{i=1}^\infty$. 
Moreover we can set $A^1$ to zero wherever the conclusion
of Lebesgue differentiation theorem \eqref{lebdiff} fails.
We claim that $\supp A^1\cap F(\supp A^2)=\emptyset$.

Suppose that 
is not the case. Then there are points $a_2\in \supp A_2$
and $a_1\in \supp A_1$ such that $\abs{a_2-a_1}\in \supp \sigma_2$.
Pick a small enough $\delta$ so that
$Q(0,\delta)\cap \supp \sigma_2=\emptyset$.
Then $\bigl(Q(0,\delta)\cap A^1_i\bigr)\cup A^2_i$ is
$\sigma_2$-avoiding for every $i$. Since the
conclusion of Lebesgue differentiation theorem
holds for every point of 
$\bigl(Q(0,\delta)\cap \supp A^1\bigr)\cup \supp A^2$ by
the argument of theorem \ref{weakstarlim} the
set $\bigl(Q(0,\delta)\cap \supp A^1\bigr)\cup \supp A^2$
is $\sigma_2$-avoiding. This proves the claim.

Furthermore, $A^1(x)\leq m(P_1)$ for all $x$. Indeed, suppose
$A^1(x_0)\geq m(P_1)(1+\epsilon)$. Since $A^1$ satisfies the conclusion of 
Lebesgue differentiation theorem \eqref{lebdiff} at $x_0$,
we can choose $\delta$ small enough so that $\delta^{-d}\int_{Q(x_0,\delta)}A^1(x)\,dx
\geq m(P_1)(1+\epsilon/2)$. Then for all sufficiently large
$i$ by lemma~\ref{ratefact} we 
have 
\begin{equation*}
d_{Q(x_0,\delta)}(A_i^1)\geq m(P_1)\left(1+\frac{\epsilon}{3}\right)
\geq m_{Q(x_0,\delta)}(\tfrac{1}{t_i}\cdot P_1)
\left(1+\frac{\epsilon}{3}\right)\Big/\left(1+\frac{\diam P_1}{t_i\delta}\right)^d
> m_{Q(x_0,\delta)}(P_1)
\end{equation*}
which is in contradiction with the fact that $A^1_i$ is $\sigma_1$-avoiding.

Let $\alpha=d_{Q(0,R)}\bigl(F(\supp A^2)\bigr)$.
Therefore by the lemma above 
\begin{align*}
\lim_{i\to\infty} d_{Q(0,R)}(A_i)&=\frac{1}{\meas{Q(0,R)}}
\int \bigl(A^1(x)+A^2(x)\bigr)\,dx\\
&\leq \alpha m_{Q(0,R)}(P_2)\left(1+\frac{\diam P_2}{R}\right)^d+
(1+\diam P_2/R)^d-1+(1-\alpha)m(P_1)\\
&\leq \max\bigl(m(P_1),m(P_2)\bigr)+\epsilon.
\end{align*}
Since $\epsilon$ was arbitrary, the proof is complete.
\end{proof}

\begin{proof}[Proof of theorem \ref{computthm}]
Note that the proofs of zooming-out lemma and supersaturation lemma are 
effective: 
the dependencies between all the constants are effectively computable.

For integer $k$ partition $Q(0,1)$ into $k^d$ cubes in the natural way.
Say a set $A$ is $k$-granular if $A$ is union of some of these cubes.
Let $\mathcal{G}_k$ be the collection of $k$-granular sets.
Let $\mathcal{P}_k=\{A+\Z^d : Z\in \mathcal{G}_k\}$.
The following simple algorithm outputs $m(D)$ within
absolute error $8\epsilon$.
\renewcommand*{\theenumi}{\arabic{enumi}}
\begin{enumerate}
\item If $\epsilon>1/10$, set $\epsilon=1/10$.
\item Make the following assignments:
\begin{align*}
r&=\min D,\\
\tilde{m}&=\left(\frac{rd^{-1/2}}{rd^{-1/2}+\diam D}\right)^d,\\
\epsilon_2&=\epsilon^{2d},\\
\epsilon_1&=\tilde{m}\epsilon^{2d}.
\end{align*}
\item Set $\delta_0$ and $R_0$ to the values
whose existence is asserted by the weak supersaturation
lemma with $\epsilon_1$ and $\epsilon_2$ as above.
\item Make the following assignments
\begin{align*}
R&=\max\left(R_0,\frac{d\diam D}{\epsilon^{2d}}\right),\\
k&=\max\left(\lceil 1/\epsilon\rceil,\lceil 1/\delta_0\rceil \right).
\end{align*}
\item Let $P$ be the property of being $(1/R)\cdot D$-avoiding.
By checking each set in $\mathcal{P}_{k^2}$ compute 
\begin{equation*}
m'=\max_{\substack{A\in \mathcal{P}_{k^2}\\P(A)=1}}d(A).
\end{equation*}
\item Output $m'$.
\end{enumerate}
The first step of the algorithm allows us
to assume that $\epsilon\leq 1/10$ in our analysis.
Note that since $Q(0,rd^{-1/2})+(rd^{-1/2}+\diam D)\Z^d$
is $D$-avoiding, we have $\tilde{m}\leq m(D)$. Clearly,
$m'\leq m(t\cdot D)=m(D)$. We will show that $m'\geq m(D)(1-8\epsilon)$.

By theorems \ref{locoptthm} and \ref{density_indep}
there is a $D$-avoiding set of density~$m(D)$.
Thus, by the proof of 
lemma~\ref{ratefact} 
there is a periodic $D$-avoiding 
set $A$ with the period $R$ and density 
$d(A)\geq m(D)/\left(1+\epsilon^{2d}/d\right)^d\geq 
m(D)(1-\epsilon^{2d})$. 
Let $A'=\Zm_{1/k}(1-\epsilon^d)A$.
If $d(A')\leq m(D)(1-3\epsilon^d)$ then
\begin{equation*}
m(D)(1-\epsilon^{2d})\leq d(A)\leq \epsilon^d d(A')+
d(\Zm_{1/k}(\epsilon_1)A)(1-\epsilon^d)+
\epsilon_1
\end{equation*}
implies that 
\begin{equation*}
d(\Zm_{1/k}(\epsilon_1)A)\geq \frac{m(D)(1-\epsilon^d+\epsilon^{2d})}{1-\epsilon}
\geq m(D)(1+\epsilon^{2d})
\end{equation*}
and weak supersaturation lemma tells us that $A$ is
not $D$-avoiding. Thus $d(A')\geq m(D)(1-3\epsilon^d)$.

Consider the set $(A'+x)\cap (R/k)\Z^d$ for a vector $x$
chosen uniformly at random from $Q(0,R/k)$.
By averaging there is a choice of $x$ for which 
$\abs{(A'+x)\cap (R/k)\Z^d}\geq d(A')k^d$. Let $x_0$ be such a choice.
Set $A''=\bigl((A'+x_0)\cap (R/k)\Z^d\bigr)+
Q\bigl(0,(1-3^{1/d}\epsilon)R/k\bigr)$. The
set $A''$ is $D$-avoiding. 
Indeed, suppose for some $x,y\in A''$ we have
$\abs{x-y}\in D$. Then
$Q(x,3^{1/d}\epsilon R/k)$ is contained
in a cube of side length $1/k$ on which $A+x_0$ has density
at least $1-\epsilon^d$. Let $\tau=d_{Q(x,3^{1/d}\epsilon R/k)}(A+x_0)$.
Then $1-\epsilon^d\leq 1-(1-\tau)(3^{1/d}\epsilon)^d$
implying $\tau\geq 2/3$.
Similarly, for $d_{Q(y,3^{1/d}\epsilon R/k)}(A+x_0)\geq 2/3$.
Therefore $A$ is not $D$-avoiding. This contradiction
proves that $A''$ is $D$-avoiding.

The set $(1/R)\cdot A''$ is $(1/R)\cdot D$-avoiding
and periodic with the fundamental region $Q(0,1)$.
It is also a union of cubes
of side length $(1-3^{1/d}\epsilon)/k$.
Each such cube contains
$k^2$-granular set 
of measure at least $[(1-3^{1/d}\epsilon)/k-2/k^2]^d$. 
Therefore $A''$ contains
a $k^2$-granular set of density
at least $d(A'')\left(1-\frac{2}{k(1-3^{1/d}\epsilon)}\right)
\geq d(A'')(1-3/k)$ since $\epsilon\leq 1/10$.
Thus there is a $k^2$-granular set
of density at least 
$m(D)(1-3\epsilon)(1-3^{1/d}\epsilon)(1-3\epsilon^d)\geq 
m(D)(1-8\epsilon)$.
\end{proof}

\section{Concluding remarks}
Let $G$ be a finite graph, and suppose
that for every edge $e\in E(G)$ there is
an admissible measure $\sigma_e\in\mathcal{M}(\R^d)$.
Then we say that a copy of the graph $G$ occurs in a set
$A\subset \R^d$ if there is a map $f\colon V(G)\to A$ such
that for every edge $xy\in E(G)$ we have
$f(x)-f(y)\in\supp \sigma_{xy}$. The theorems
\ref{distsatur} and \ref{ORsatur} show that
if $G$ is either a single edge or a path of
length~$2$, then the property of avoiding
$G$ is supersaturable. The proof of theorem
\ref{ORsatur} can be easily modified to
the case when $G$ is a star. I conjecture
that the property of avoiding $G$ is supersaturable
whenever $G$ is a tree. 

An example of Bourgain
\cite{cite:bourgain} shows that the property
of avoiding a triangle $K_3$ fails to
be supersaturable. However, in his example
the points of the triangle are forced to
lie on the same line. Perhaps with an
appropriate non-degeneracy condition the
property of avoiding $K_3$ is supersaturable.

Suppose $G_1$ and $G_2$ are two graphs as above.
Let $r_1,r_2$ be two distinguished vertices in
$G_1$ and $G_2$ respectively. Then
$G_1\OR G_2$ is a graph which is obtained
by identifying $r_1$ and $r_2$ in the
disjoint union of $G_1$ and $G_2$. I believe
that in the case when $G_1$ and $G_2$ are trees,
the generalization of theorem \ref{ORthm} holds:
$m(G_1\OR t\cdot G_2)\to \max\bigl(m(G_1),m(G_2)\bigr)$
as $t\to\infty$.

Further problems on configurations in sets of positive
density and the survey of known results
can be found in \cite[\S6.3]{cite:resprobgeom}.

The theorem \ref{smpand} implies that the measurable
chromatic number $\chi_{\R^d}^m(D)$ grows
exponentially in $\abs{D}$ provided the elements of $D$ grow
fast enough. It is very likely that the usual
chromatic number does not share this behavior.
I conjecture that for any dimension $d$ if the elements of $D$ are algebraically
independent over $\Q$, then 
the chromatic number $\chi_{\R^d}(D)$ is
bounded independently of what $D$ actually is. The conjecture
is easily seen to be true when $d=1$ because the finite subgraphs of $G_{\R^1}(D)$
are subgraphs of the $\abs{D}$-dimensional rectangular grid \raisebox{-0.15em}{\includegraphics[scale=0.05]{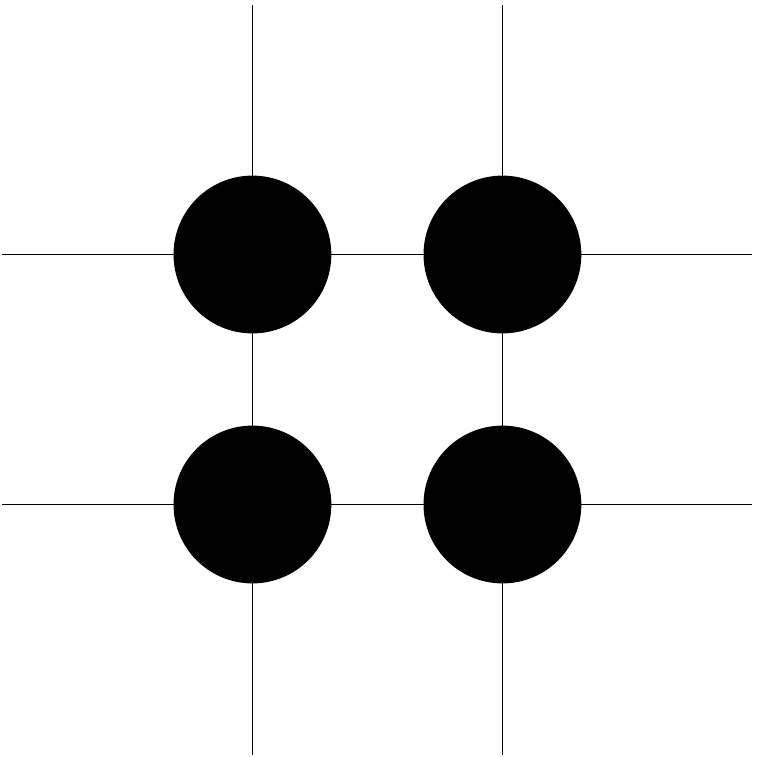}}. 
The \emph{clique number} of a graph $G$, denoted $\omega(G)$, is the
number of vertices in the largest complete subgraph of~$G$.
For $d\geq 2$ the only result that I can prove is
\begin{theorem}
There is a function $f(d)$ such that
if the elements of $D$ are algebraically independent over~$\Q$, then 
$\omega\bigl(G_{\R^d}(D)\bigr)<f(d)$.
\end{theorem}
\begin{proof}
Denote by $K_n$ the complete graph on $n$~vertices.
Suppose $K_n$ is a subgraph of $G_{\R^d}(D)$. Then
let $X=\{x_1,x_2,\dotsc,x_n\}\subset \R^d$ be the vertices of 
this complete subgraph. 
Let $A=(a_{i,j})_{i,j=1}^n$ be an $n\times n$ matrix whose
entries are $a_{i,j}=\dist(x_i,x_j)^2=\langle x_i-x_j,x_i-x_j\rangle
=\langle x_i,x_i\rangle+\langle x_j,x_j\rangle-2\langle x_i,x_j\rangle$.
The matrix $B=(b_{i,j})_{i,j=1}^n$ with
$b_{i,j}=\langle x_i,x_i\rangle$ has rank~$1$. The matrix
$C=(c_{i,j})_{i,j=1}^n$ with $c_{i,j}=\langle x_i,x_j\rangle$
has rank at most~$d$. Thus the rank of $A=B+B^{t}-2C$ is at most~$d+2$.

Consider any subset $X'\subset X$ of $d+3$ elements.
Let $A'$ be the corresponding $(d+3)\times (d+3)$ submatrix 
of~$A$.
Let $r_1,\dotsc,r_k$ be the non-zero elements that 
occur in $A'$. Since $r_1,\dotsc,r_k$ are
squares of algebraically independent numbers, they
themselves are algebraically independent.
Since $A'$ is not of the full rank, $\det A'=0$.
Since the determinant is a polynomial function with rational
coefficient in entries of $A'$, it follows that
$\det A'=0$ whenever $\{r_1,\dotsc,r_k\}$
is replaced by any set of $k$ algebraically
independent numbers. Therefore, $\det A'$
is zero as a polynomial in $r_1,\dotsc,r_k$.
Since the matrix $A'$ is a symmetric matrix,
each $r_i$ occurs at least twice. If 
each $r_i$ occurred exactly twice, then $\det A'(r_1,\dotsc,r_k)$,
being the determinant of the general symmetric matrix
with the zeros on the diagonal, would not be the zero polynomial. 
Thus, in every
set of $d+3$ points at least one distance occurs
twice.

Color the edge $x_ix_j$ of the complete graph on $X$
by the distance between $x_i$ and~$x_j$. The above
asserts that there is no $K_{d+3}$ subgraph whose edges 
all  colored differently.
On the other hand, since the simplex on $d+2$ vertices
does not embed isometrically in $\R^d$,
there is no monochromatically colored $K_{d+2}$ subgraph.
By the canonical Ramsey theorem \cite{cite:canonical_ramsey}
if $n$ is large enough, then there is 
a $Y= \{y_1,\dotsc,y_{d+4}\}\subset X$ 
such that the color of an edge $y_iy_j$ for $i<j$ depends only
on~$i$. Let $t_i=\dist(y_i,y_{i+1})^2$. The $(d+4)\times
(d+4)$ matrix corresponding to $Y$ is 
$M=(m_{i,j})_{i,j=1}^{d+4}$ where 
\begin{equation*}
m_{i,j}=
\begin{cases}
t_i,&\text{if }i<j,\\
t_j,&\text{if }i>j,\\
0,&\text{if }i=j.
\end{cases}
\end{equation*}
The matrix $M$ is of rank at least $d+3$. Indeed,
let $r_i$ be the $i$'th column of $M$. Then for every
$i=1,\dotsc,d+3$
the first $i-1$ coordinates of $r_{i+1}-r_i$ are zero,
and $i$'th coordinate is non-zero. Thus the vectors
$r_{i+1}-r_i$ span a vector space of dimension $d+3$ implying
that $M$ is of rank at least $d+3$. Since
$M$ is a submatrix of $A$, which is of rank at most $d+2$,
we reached a contradiction.
\end{proof}

\textbf{Acknowledgements.} I would like to thank Uri Andrews and 
Pablo Candela-Pokorna for discussions that inspired this work.
I am very grateful to Josef Cibulka and Jan Kyn\v{c}l
for a careful reading of an earlier version of the paper and many 
useful suggestions. 

\bibliographystyle{alpha}
\bibliography{forbdist}

\end{document}